\documentclass[8pt]{article}
\usepackage{amsfonts}
\usepackage{mathrsfs}
\usepackage{amsmath}
\usepackage{dsfont}
\usepackage{color}
\usepackage{graphicx,epstopdf}
\usepackage{amssymb}
\usepackage{theorem}
\usepackage{graphicx}
\usepackage[dvips]{epsfig}
\usepackage{latexsym}
\usepackage{exscale}
\usepackage[latin1]{inputenc}
\newtheorem{theorem}{Theorem}[section]

\newtheorem{corollary}[theorem]{Corollary}

\newtheorem{definition}[theorem]{Definition}
\newtheorem{example}[theorem]{Example}

\newtheorem{lemma}[theorem]{Lemma}

\newtheorem{proposition}[theorem]{Proposition}
\newtheorem{remark}[theorem]{Remark}

\def\ud{\, \mathrm{d}}
\numberwithin{equation}{section}

\begin{document}

\title{Fractional Hida Malliavin Derivatives and Series Representations of Fractional Conditional Expectations}
\date{}
\maketitle
\centerline{\bf {Sixian Jin\footnote{Institute of Mathematical Sciences, Claremont Graduate University, Sixian.Jin@cgu.edu}, Qidi Peng\footnote{Institute of Mathematical Sciences, Claremont Graduate University, Qidi.Peng@cgu.edu}, Henry Schellhorn\footnote{Corresponding author, Institute of Mathematical Sciences, Claremont Graduate University, Henry.Schellhorn@cgu.edu} }}

\centerline{}

\begin{abstract}
We represent fractional conditional expectations of a
functional of fractional Brownian motion as a convergent series in $L^2(\mathbb P^H)$ space. When the target random variable is some function of a discrete trajectory of fractional Brownian motion, we obtain a backward Taylor series representation;  when the target functional is generated by a continuous fractional filtration, the series representation is obtained by applying a "frozen path" operator and an exponential operator to the functional. Three examples are provided to show that our representation gives useful series expansions of ordinary expectations of target random variables.
\end{abstract}

\section{Introduction}

Fractional Brownian motion (fBm), nowadays considered as one of the most natural
extensions of the classical Brownian motion (Bm), was first introduced by
Kolmogorov in 1940 \cite{Kolmogorov}. It is popularized by Mandelbrot and Van
Ness \cite{Mandelbrot} in 1968, by being introduced to mathematical
finance, to take into account the effect of exogenous arrivals. Today, the mathematical analysis and applications of fBms have advanced enormously in a variety of fields and at all levels, such as stock returns modeling, fractal geometry, signal processing,
geology, geography and statistical biology, etc. FBm is almost applicable in every domain where one can apply Bm.  However, unlike Bm, some new mathematical issues arise when replacing
Bm by such an extension, which is mainly due to the complex
covariance structure of its increments. For example, fBm has correlated increments for $H\neq 1/2$. This results that evaluating conditional expectations of functions or functionals of
fBm $\{B_{t}^{H}\}_{t\in \mathbb R}$ is a notoriously difficult problem. Grippenberg and Norros \cite{Gripenberg} provided a
technical and difficult approach to calculate the conditional mean of fBm. Fink et al. \cite{Fink} also addressed this problem when studying the price of a
zero-coupon bond in a fractional bond market. Since fBm is generally not a Markov process, both authors restricted themselves to calculate conditional expectations given the current value of $B_{t}^{H}$, and not given the whole path of $B_{t}^{H}$ preceding $t$.

This paper presents a different and original way to evaluate expectations of
functionals of fBm. It is based on the Malliavin
calculus with respect to fBm as presented in Biagini et al. \cite%
{Biagini}, Chapter 2 and Chapter 3. Since the 70's, Malliavin calculus plays an important role in analysis on the Wiener space as well as in the study on stochastic differential equations. The main advantage of Malliavin
calculus is that, it allows to give sufficient conditions for the distribution of a random
variable to have a smooth (differentiable) density with respect to
Lebesgue measure and to give bounds for this density and its derivatives. With this technique, one of the main fruits from Malliavin calculus is the fractional conditional expectation and its representation by Clark-Hausmann-Ocone formula (see e.g. \cite{Biagini}). By using the so-called fractional Clark-Hausmann-Ocone formula, our first main result represents the fractional conditional expectations of functions of fBm's discrete trajectory as a convergent series in $L^2(\mathbb P^H)$. Our second main result is more general under a different sufficient condition for the convergence. It is obtained from the fact that, the fBm has the  "martingale" property under fractional conditional expectation. This property leads to an exponential expansion of the fractional conditional expectations. The latter result partially extends the work of Schellhorn and Jin \cite{Schellhorn}, who proved this representation for conditional expectations of a functional of Bm. It is worth noting that,
the fBm can be divided into three very different
classes according to the values of $H>1/2$, $H=1/2$ and $H<1/2$. When $H>1/2$, the corresponding fBm is persistent, which means that its increments are positively correlated (the increase of the increments is
likely to be followed by another increase). A huge number of phenomena can
be modeled in terms of this class of processes, such as the level of the optimum dam
sizing, the logarithm of the stock return and
financial turbulence \cite{Mandelbrot}. As in \cite{Biagini}, Chapter 3, we only consider the case where $%
H>1/2$ in this work. We remark that a similar study can be done
for classes of fBm with $H<1/2$ for the future.

The structure of this paper is as follows. In Section 2, we
introduce some definitions, notations and known results on Malliavin
calculus related to fBm, which are needed for the next sections. In Section 3, we present our two main results, which are i) a
generalization to fBm of the backward Taylor
expansion obtained in \cite{Schellhorn}, and ii) the exponential formula
itself. We note that both series representations can be used both for
numerical applications and for solving some algebraic problems. In section 4, we show 3 applications of the exponential formula, respectively to the fractional Merton model of interest rates, to a special case of the fractional Cox-Ingersoll-Ross model of interest rates, and to the characteristic function of geometric fBm.
\section{Preliminaries}
\subsection{Fractional Brownian Motion}
A real-valued standard fBm can be defined independently and equivalently using a moving average representation \cite{Mandelbrot} and a harmonizable representation \cite{Samorodnitsky}.  In fact, up to a multiplicative scaling factor, fBm is the unique Gaussian, self-similar, with stationary increments process. Therefore, a standard fBm can be defined from the uniqueness of its covariance structure:
\begin{definition}
\label{definition 2.1.}
   A standard fBm $\{B_t^H\}_{t\in\mathbb R}$ with Hurst index $H\in(0,1)$ is the unique centered Gaussian process with almost surely continuous non-differentiable sample path and with covariance function: for any $s,t\in\mathbb R$,
\begin{equation*}
E\left[ B_{s}^{H}B_{t}^{H}\right] =\frac{1}{2}\Big(
|s|^{2H}+|t|^{2H}-\left\vert t-s\right\vert ^{2H}\Big).
\end{equation*}
\end{definition}
Without any loss of generality, we restrict the fBm to nonnegative-time process and denote the corresponding probability space by $\left(
\Omega ,(\mathcal{F}_{t}^{H})_{t\in\mathbb R_+},\mathbb{P}^{H}\right) $, where $(\mathcal{F}%
_{t}^{H})_{t\in\mathbb R_+}$ is the natural filtration generated by fractional Wiener chaos (see e.g. \cite{Biagini}, Page 49) and $\mathbb{P}^{H}$ is the corresponding probability
measure. Again, as an assumption, we let $H>1/2$ in the remaining of the paper.
\subsection{Fractional Hida Malliavin Derivative}
Let $L^2(\mathbb P^H):=L^2(\Omega,(\mathcal F_t^H)_{t\in\mathbb R_+},\mathbb P^H)$ denote a Hilbert space of random variables $F$ equipped with the norm
 $$
 \|F\|_{L^2(\mathbb P^H)}:=\sqrt{E|F|^2}<+\infty.
 $$
 We indicate by $\mathcal{S}\left( \mathbb{R_+}\right) $ the Schwartz space of
rapidly decreasing smooth functions on $\mathbb{R_+}$. More precisely,
$$
\mathcal{S}( \mathbb{R_+}):=\Big\{f\in C^{\infty}(\mathbb R_+):~\mbox{for any}~ \alpha,\beta\in\mathbb N,~\sup_{x\in\mathbb R_+}\Big|x^\alpha\frac{\ud^\beta f(x)}{\ud x^\beta}\Big|<+\infty\Big\},
$$
where $C^{\infty}(\mathbb R_+)$ denotes a space of continuously infinitely differentiable real-valued functions.
 We equip $\mathcal{S}\left( \mathbb{R_+}\right) $ with the following
inner product: for any $f,g\in \mathcal{S}\left( \mathbb{R_+}\right)$,
\begin{equation*}
\left\langle f,g\right\rangle _{H}:= \int_{\mathbb{R_+}}\int_{\mathbb{R_+%
}}f\left( s\right) g\left( t\right) \varphi_H \left( s,t\right) \ud s\ud t,
\end{equation*}
where  $\varphi_H \left( s,t\right)
:=H\left( 2H-1\right) \left\vert s-t\right\vert ^{2H-2}$ for any $s,t\in\mathbb R_+$.

Denote by $%
L_{\varphi_H }^{2}\left( \mathbb R_+\right) $ the completion of $\mathcal{S}\left(
\mathbb R_+\right)$ under the norm $\|f\|_{H,\mathbb R_+}:=\sqrt{\langle f,f\rangle _{H}}$. This is a separable Hilbert space of deterministic functions.
 Note that fBm is an isonormal process and, if two functions $f,g\in L_{\varphi_H }^{2}\left( \mathbb R_+\right)$, the stochastic integrals with respect to fBm $\int_{\mathbb R_+}f(s)\ud B_s^H$ and $\int_{\mathbb R_+} g(s)\ud B_s^H$ are well defined, zero mean, Gaussian random variables with  covariance
$$
Cov\Big(\int_{\mathbb R_+}f(s)\ud B_s^H,\int_{\mathbb R_+}g(s)\ud B_s^H\Big)=\left\langle f,g\right\rangle _{H}.
$$
Let $\mathcal{S}^{\prime }\left( \mathbb{R_+}\right) $
 denote the dual space of $\mathcal{S}\left( \mathbb{R_+}\right)$. $\mathcal{S}^{\prime }\left( \mathbb{R_+}\right)$ is also called the space of tempered distributions on $\mathbb R_+$ (see e.g. \cite{Pipiras} or \cite{Biagini}). Note that by the Bochner-Milos theorem (see e.g. \cite{Biagini,Holden,Kuo}), the probability measure $\mathbb P^H$ is the one which allows $$
 B^H_t:=<B^H,\chi_{[0,t]}>:=\int_0^t1\ud B^H_u$$ to be an element in $L^2(\mathbb P^H)$.

For the purpose of presenting an element $F\in L^2(\mathbb P^H)$ using convergent Taylor series, we first introduce the notion of fractional Hida Malliavin derivative (see Definition 3.3.1 in \cite{Biagini}).
\begin{definition}[Fractional Hida Malliavin Derivative]
\label{Definition 2.3.} Given an operator $G:\mathcal{S}^{\prime }\left( \mathbb{R_+}
\right) \rightarrow \mathbb R $ and some $\gamma\in \mathcal{S}^{\prime }\left( \mathbb{R_+}\right)$. $G$ is said to have a directional derivative in the direction $\gamma$ if, for all $\nu
\in \mathcal{S}^{\prime }\left( \mathbb{R_+}\right)$, there exists an element $X_{\nu,\gamma}$ in the  fractional Hida distribution space $\left( \mathcal{S}\right)_{H}^{\ast }$ (see  Definition 3.1.10 in \cite{Biagini}), such that
$$
\lim_{\epsilon\rightarrow0}\frac{G\left( \nu +\epsilon \gamma \right) -G\left(
\nu \right) }{\epsilon }=X_{\nu,\gamma},~\mbox{in}~\left( \mathcal{S}\right)_{H}^{\ast }.
$$
We say $G$ is fractional Hida Malliavin differentiable if there exists a map $\Psi :\mathbb{%
R_+}\times \mathcal{S}'(\mathbb{R_+})\rightarrow \left( \mathcal{S}\right)_{H}^{\ast }$ such that for all $\nu
\in \mathcal{S}^{\prime }\left( \mathbb{R_+}\right)$,
$
\Psi \left( \cdot,\nu \right)
\gamma \left( \cdot\right)
$
 is $\left( \mathcal{S}\right)_{H}^{\ast }$-integrable and
\begin{equation*}
X_{\nu,\gamma} =\int_{\mathbb{R_+}}\Psi \left( t,\nu
\right) \gamma \left( t\right) \ud t
\end{equation*}%
for all $ \gamma \in L^2(\mathbb R_+)\subset\mathcal S'\left( \mathbb{R_+}\right)$. Then we set, for all $t\in\mathbb R_+$,
\begin{equation*}
D_{t}^{H}G\left( \nu \right) := \Psi \left( t,\nu \right)
\end{equation*}%
and we call $D_{t}^{H}G\left( \nu \right) $ the fractional Hida Malliavin derivative with order $H$ of $G$ on $\nu$ at  $t$.
\end{definition}

Remark that the fractional Hida Malliavin derivative with respect to fBm extends the classical one (see \cite{Nualart}) with respect to Bm. Also remark that, since $L^2(\mathbb P^H)\subset \left( \mathcal{S}\right)_{H}^{\ast }$ (see Definition 3.1.10 in \cite{Biagini}), we will mainly focus on the fractional Hida Malliavin derivative defined on $L^2(\mathbb P^H)$ in the sequel.

It is useful to note that the fractional Hida Malliavin derivative possesses some nice properties similar to the classical derivatives. For example, the chain rule is still valid: if $G(\nu)=H\left( B_{t_1}^{H},...,B_{t_n}^{H}\right)$, with $H:~\mathbb R^n\rightarrow\mathbb R$ being some deterministic differentiable function, then for $t\ge0$,
\begin{equation}
\label{multiD}
D_t^{H}G(\nu)=\sum\limits_{i=1}^{n}\frac{\partial H}{\partial x_{i}}\left(
B_{t_1}^{H},...,B_{t_n}^{H} \right) \chi _{\left[
0,t_{i}\right] }(t).
\end{equation}
Now we introduce some other important properties of fractional Hida Malliavin derivative that we will need to construct the convergent series of fractional conditional expectations. The most interesting one is fractional Clark-Hausmann-Ocone formula, which generalizes the classical Clark-Hausmann-Ocone formula. The following statements as well as notations are helpful to the introduction of this formula.
\begin{theorem} [It\^o decomposition, see \cite{Biagini}, Page 82]
\label{Theorem 2.4.}
Let $F\in L^2(\mathbb P^H)$, then $F$ has the following representation via fractional Wick It\^o Skorohod integral (FWISI): there exists a sequence of deterministic functions $f_{n}\in \hat{L}_{\varphi_H}^{2}\left( \mathbb{R}%
_{+}^{n}\right) $ such that
\begin{equation*}
F=\sum\limits_{n=0}^{+\infty }I_{n}\left( f_{n}\right)~\mbox{in $L^2(\mathbb P^H)$ with}~
\left\Vert F\right\Vert _{L^{2}\left( \mathbb{P}^{H}\right)
}^{2}=\sum\limits_{n=0}^{+\infty }n!\left\Vert f_{n}\right\Vert _{H,\mathbb R_+^n }^{2},
\end{equation*}%
where
\begin{description}
\item[$-$] the sequence of multiple FWISI $(I_n(f_n))_{n\in\mathbb N}$ is defined in \cite{Biagini}, Page 81;
\item[$-$] $\hat{L}_{\varphi_H}^{2}\left(
\mathbb{R}_{+}^{n}\right) $ denotes the subspace of symmetric functions in $%
L_{\varphi_H}^{2}\left( \mathbb{R}_{+}^{n}\right) $;
\item[$-$] $\left\Vert f_{n}\right\Vert _{H,\mathbb R_+^n}$ is the norm defined by
$$\left\Vert f_{n}\right\Vert _{H,\mathbb R_+^n}^2:=n!\int_{\mathbb R_+^{2n}}f_n(u_1,\ldots,u_n)f_n(v_1,\ldots,v_n)\prod_{i=1}^n\varphi_H(u_i,v_i)(\ud u)^{\otimes n}(\ud v)^{\otimes n},$$
with
\begin{equation}
\label{du}
(\ud u)^{\otimes n}:=\ud u_n\ud u_{n-1}\ldots\ud u_1.
\end{equation}
\end{description}
\end{theorem}
\bigskip
\begin{definition}
We define $h_n$, the Hermite polynomial of degree $n\ge0$, by $h_0\equiv1 $ and for $n\ge1$,
\begin{equation}
\label{Hermite}
h_n(x)=(-1)^n\exp\left(\frac{x^2}{2}\right)\frac{\ud^n}{\ud x^n}\exp\left(-\frac{x^2}{2}\right),~\mbox{for all $x\in\mathbb R$}.
\end{equation}
Remark that the sequence $\{h_n\}_{n\ge0}$ is equivalently defined as solutions of the following equation, valid for all $t,x\in\mathbb R$:
\begin{equation}
\label{defnHermite}
\exp\left(tx-\frac{t^2}{2}\right)=\sum_{n=0}^{+\infty}\frac{t^n}{n!}h_n(x).
\end{equation}
\end{definition}

\begin{proposition}
\label{Proposition 2.5.}
Let $f\in\hat{L}_{\varphi_H}^{2}\left( \mathbb{R}%
_{+}^{n}\right) $, then the multiple FWISI $I_n(f)$ exists and is given as
\begin{equation*}
I_{n}\left( f\right) =n!\int_{0\leq s_{1}\leq \ldots\leq s_{n}<+\infty}f\left(
s_{1},\ldots,s_{n}\right) \ud \left( B_{s}^{H}\right) ^{\otimes n},
\end{equation*}
where we denote by
\begin{equation}
\label{dBs}
 \ud \left(B_{s}^{H}\right) ^{\otimes n}:=\ud B_{s_n}^H\ud B_{s_{n-1}}^H\ldots\ud B_{s_1}^H.
 \end{equation}
In particular, if there exists $g\in L_{\varphi_H}^2(\mathbb R_+)$ such that $f\left(
s_{1},\ldots,s_{n}\right) =$$g(s_{1})\ldots g(s_{n})$ for all $(s_1,\ldots,s_n)\in\mathbb R_+^n$, then
\begin{equation*}
I_{n}\left( f\right) =\left\Vert g\right\Vert _{H,\mathbb R_+}^{n}h_{n}\left( \frac{%
\int_{0}^{+\infty}g\left( s\right) \ud B_{s}^{H}}{\left\Vert g\right\Vert _{H,\mathbb R_+}}%
\right),
\end{equation*}
where the norm $\|g\|_{H,\mathbb R_+}$ equips the space $L_{\varphi_H}^2(\mathbb R_+)$.
\end{proposition}
As a special case, when taking $f=\chi_{[t,T]^n}$ with $0\leq t<T$ in Proposition \ref{Proposition 2.5.}, we obtain
\begin{equation}
\label{eq}
\int_{t\leq s_{1}\leq\ldots\leq s_{n}\leq T}1 \ud \left(B_{s}^{H}\right) ^{\otimes n}=%
\frac{\left( T-t\right) ^{nH}}{n!}h_{n}\left( \frac{B_{T}^{H}-B_{t}^{H}}{%
\left( T-t\right) ^{H}}\right).
\end{equation}
From now on, the FWISI of a continuous-time stochastic process $\{X(s)\}_{s\ge0}$ over any time interval $[a,b]$, is denoted by $\int_a^bX(s)\ud B_s^H$.
\begin{definition}[Fractional Conditional Expectation]
\label{Definition 2.6.}
Let $F\in L^2(\mathbb P^H)$ be represented as (such an expansion exists, due to Definition 3.10.1 in \cite{Biagini})
$$
F=\sum\limits_{n=0}^{+\infty }\int_{\mathbb{R}_{+}^{n}}g_{n}%
\left( s_1,\ldots,s_n\right) \ud\left( B_s^{H}\right) ^{\otimes n},
$$
with some sequence of functions $g_n\in\hat{L}_{\varphi_H}^{2}\left(
\mathbb{R}_{+}^{n}\right)$.  Then for $t\ge0$, we define the fractional
conditional expectation of $F$ with respect to $\mathcal{F}_{t}^{H}$ by%
\begin{equation}
\label{tildeE}
\tilde{E}\left[ F|\mathcal{F}_{t}^{H}\right] :=
\sum\limits_{n=0}^{+\infty }\int_{[0,t]^n}g_{n}\left( s_1,\ldots,s_n\right)\ud\left( B_s^{H}\right) ^{\otimes
n}.
\end{equation}
\end{definition}
\textbf{Remarks:}
Though different from conditional expectation, fractional conditional expectation has some properties, which are similar to those of classical conditional expectation:
\begin{enumerate}
\item For all $t\ge0$, $\tilde{E}\left[ F|\mathcal{F}_{t}^{H}\right]$ is $\mathcal F_t^H$-measurable. Vice versa, if $F$ is $\mathcal F_t^H$-measurable,
$\tilde{E}\left[ F|\mathcal{F}_{t}^{H}\right]=F$.
    \item For any $s, t\ge0$,
\begin{equation}
\label{interchangetildeE}
\tilde{E}\left[\tilde{E}\left[ F|\mathcal{F}_{s}^{H}\right]|\mathcal F_t^H\right]=\tilde{E}\left[ F|\mathcal{F}_{\min\{s,t\}}^{H}\right].
\end{equation}
\item  From (\ref{tildeE}), for all $t\ge0$, $\tilde{E}\left[ F|\mathcal{F}_{t}^{H}\right]$ is $\mathcal F_t^H$-measurable. By virtue of Lemma 3.10.5 1)  in \cite{Biagini} and (2.14), (2.22) in \cite{Aase} (by taking $r=0$), one has,
\begin{equation}
\label{boundtildeE}
\left\|\tilde{E}\left[ F|\mathcal{F}_{t}^{H}\right]\right\|_{L^2(\mathbb P^H)}\leq \left\|F\right\|_{L^2(\mathbb P^H)}.
\end{equation}
This inequality yields that, if $F\in L^2(\mathbb P^H)$, then $\tilde{E}\left[ F|\mathcal{F}_{t}^{H}\right]\in L^2(\mathbb P^H)$.
\item Under the transform $\tilde E[\cdot|\mathcal F_t^H]$, one recovers a property very similar to that of classical martingale. Although fBm is generally not a martingale, it is shown that, for $0\leq t\leq T$, $\tilde{E}\left[ B_{T}^{H}|\mathcal{F}_{t}^{H}\right] =B_{t}^{H}$ $\mathbb P^H$-a.s..
\end{enumerate}
The following theorem, given in \cite{Biagini}, extends the Clark-Hausmann-Ocone Formula from Bm to fBm.
\begin{theorem} [Fractional Clark-Hausmann-Ocone Formula]
\label{Theorem 2.7.}
Fix $T>0$, let the random variable $F\in L^{2}(\mathbb{P}^{H})$ be $\mathcal{F}_{T}^{H}$-
measurable and Hida Malliavin differentiable, then $
\tilde{E}\left[ D_{t}^{H}F|\mathcal{F}_{t}^{H}\right] \in L^{2}\left(
\mathbb{P}^{H}\right) \mbox{ for all } t\in \left[ 0,T\right]$ and%
\begin{equation*}
F=E\left[ F\right] +\int_{0}^{T}\tilde{E}\left[ D_{t}^{H}F|\mathcal{F}%
_{t}^{H}\right] \ud B_{t}^{H}.
\end{equation*}
\end{theorem}
\section{Main Results}
\subsection{Series Representation via Backward Taylor Expansion}
\begin{definition}
\label{DefinitionDF}
 Fix $T>0$. We call $\mathbb{D}_{\infty,T }^{H}$ the set of $
\mathcal{F}_{T}^H$-measurable random
variables which are infinitely fractional Hida Malliavin differentiable. Moreover, for any integer $n\ge1$,
$$
E\Big[\sup_{s_1,\ldots,s_n\in[0,T]}\big|D_{s_{n}}^{H}\ldots D_{s_{1}}^{H}F\big|\Big]^2<+\infty.
$$
To simplify notation, we denote by $D_u^{H,0}=Id$ (identity function) and by $D_u^{H,k}=\underbrace{D_{u}^{H}\circ \ldots\circ D_{u}^{H}}_{k-\mbox{tuple}}$ the $k$-th composition of the fractional Hida Malliavin derivative.
\end{definition}
\begin{definition}
\label{multiintegral}
 Let $F\in\mathbb D_{\infty,T}^H$. Assume $F=H(B_{t_1}^H,\ldots,B_{t_J}^H)$ for some integer $J\ge1$ and $0=t_0<t_1<t_2<\ldots<t_{J-1}<t_{J}=T$. $H:~\mathbb R^J\rightarrow\mathbb R$ is an infinitely differentiable deterministic function. For $j\in\{1,\ldots,J\}$ and $r\in[t_{j-1}, t_j]$,
let $\{\psi_{k}^{(r,t_j)}(F)\}_{k\in\mathbb N} $ be the sequence given as:
\begin{equation}
\label{multiS1}
\psi_{0}^{(r,t_j)}(F)=F;
\end{equation}
and for $k\ge1$, $\psi_{k}^{(r,t_j)}(F)$ equals
\begin{eqnarray}
\label{multiS2}
&&2^{-k}\sum_{\sum_{l=1}^Jq_l=k}\prod_{i=1}^J\left(\frac{\big(|t_j-t_{i-1}|^{2H}-|t_j-t_i|^{2H}+|t_i-r|^{2H}-|r-t_{i-1}|^{2H}\big)^{q_i}}{q_i!}\right) \nonumber \\
&&~~\times D_{t_1}^{H,q_1}\ldots D_{t_J}^{H,q_J}F.
\end{eqnarray}
\end{definition}
 Our first main result is given by the following theorem, where a sufficient condition is provided, such that the fractional backward Taylor expansion of the fractional conditional expectation of $F$ converges in $L^2(\mathbb P^H)$:\\
 \textbf{Assumption $(\mathcal A)$:} Let $F$  be given as in Definition \ref{multiintegral}. For some $r \in [0,T]$, assume $F$ satisfies the following:
\begin{eqnarray*}
&&\sum_{i=0}^{N}\left\|\sup_{ \sum_{j=1}^Jw_{j}=2N-i}\left|
D_{t_{1}}^{H,w_{1}}\ldots D_{t_J}^{H,w_J}F\right|\right\|_{L^2(\mathbb P^H)}{N\choose i}%
^{2}\frac{(i!)^{1/2}(N+J-1)!}{2^{N-i}(N!)^2}\\
&&~~\times\left(T^{2H}-r^{2H}+(T-r)^{2H}\right)^{N-i}(T-r)^{iH}\xrightarrow[N\rightarrow+\infty]{}0.
\end{eqnarray*}
Contrary to appearance, this condition is not difficult to check in practice. For example, any $F$ verifying
$$
 \left\|\sup_{ \sum_{j=1}^Jw_{j}=N}\left|
D_{t_{1}}^{H,w_{1}}\ldots D_{t_J}^{H,w_J}F\right|\right\|_{L^2(\mathbb P^H)}\le cN^{(1/4-\epsilon)N}
$$
for some $c>0$, some $\epsilon\in(0,1/4]$ and for all $N\in\mathbb N$ will satisfy Assumption $(\mathcal A)$. Because if so, by taking
 $$
 c_1=\frac{T^{2H}-r^{2H}+(T-r)^{2H}}{2}~\mbox{and}~c_2=T-r,
 $$the finite partial sum in Assumption $(\mathcal A)$ can be upper bounded by the following item:
\begin{eqnarray*}
&&c\sum_{i=0}^{N}{N\choose i}%
^{2}\frac{(i!)^{1/2}(N+J-1)!}{(N!)^2}(2N-i)^{(1/4-\epsilon)(2N-i)}c_1^{N-i}c_2^{i} \\
&&\le c{N\choose \lceil N/2\rceil }^{2}\frac{(N!)^{1/2}(N+J-1)!(2N)^{(1/2-2\epsilon)N}}{(N!)^{2}}\sum_{i=0}^{N}c_1^{N-i}c_2^{i}\\
&&=c\frac{(N!)^{1/2}(N+J-1)!(2N)^{(1/2-2\epsilon)N}c_1^N}{(\lceil N/2\rceil!)^2((N-\lceil N/2\rceil)!)^2}\frac{1-(c_2 c_1^{-1})^{N+1}}{1-c_2c_1^{-1}}\\
&&\thicksim \frac{(2\pi N)^{1/4}(N/e)^{N/2}\sqrt{2\pi (N+J-1)}((N+J-1)/e)^{N+J-1}}{\pi^2 N^2(N/(2e))^{2N}}(2N)^{(1/2-2\epsilon)N}c_1^N\frac{1-(c_2 c_1^{-1})^{N}}{1-c_2c_1^{-1}}\\
&&\thicksim\frac{N^{J-9/4}c_3^N}{N^{2\epsilon N}}\xrightarrow[N\rightarrow+\infty]{}0,
\end{eqnarray*}
where $\lceil\cdot\rceil$ denotes the ceiling number (the smallest integer upper bound); $c_3>0$ is some proper constant; and the last approximation is due to the Stirling's approximation.
\begin{theorem} [Fractional Backward Taylor Expansion]
\label{Theorem 3.2.} Let $F$ satisfy Assumption $(\mathcal A)$. Define
  \begin{equation}
  \label{Ir}
  I_r:=\left\{\begin{array}{ll}
  &1,~\mbox{if $r\in[0,t_1]$}\\
  &i,~\mbox{if $r\in(t_{i-1},t_i]$,~for $2\le i\le J$}
  \end{array}\right.
  \end{equation}
    and $(\tilde t_{I_r-1},\tilde t_{I_r},\ldots,\tilde t_J):=(r,t_{I_r},\ldots,t_J)$. Then the following series is convergent in $L^2(\mathbb P^H)$:
\begin{eqnarray}
\label{Ii}
&&\tilde{E}\big[ F|\mathcal{F}_{r}^{H}\big]=\sum\limits_{l=0}^{%
+\infty }\left( -1\right) ^{l}\sum_{q_{I_r}+\ldots+q_{J}=l}\sum_{i_{I_r}=0}^{q_{I_r}}\ldots\sum_{i_J=0}^{q_J}\prod\limits_{k=I_r}^{J}\frac{\left( -1\right) ^{i_k}(\tilde t_{k}-\tilde t_{k-1}) ^{\left( q_{k}-i_k\right) H}%
}{\left( q_{k}-i_k\right) !}\nonumber\\
&&\times \left(h_{q_{I_{r}}-i_{I_{r}}}\bigg( \frac{B_{\tilde t_{I_{r}}
}^{H}-B_{\tilde t_{I_{r}-1} }^{H}}{(\tilde t_{I_{r}}-\tilde t_{I_{r}-1}) ^{H}}\bigg)\psi_{i_{I_r}}^{(\tilde t_{I_r-1},\tilde t_{I_r})}\right)\circ \ldots \circ \left(h_{q_{J}-i_J}\bigg( \frac{B_{\tilde t_J
}^{H}-B_{\tilde t_{J-1} }^{H}}{(\tilde t_J-\tilde t_{J-1}) ^{H}}\bigg)\psi_{i_J}^{(\tilde t_{J-1},\tilde t_{J})} \right) \nonumber \\
&&~~(D_{t_{I_r}}^{H,q_{I_r}}\ldots D_{t_J}^{H,q_J}F).
\end{eqnarray}
In particular, when $F=H(B_T^H)$,
\begin{equation}
\label{Iii}
\tilde{E}\big[F|\mathcal{F}_{r}^{H}\big]=\sum\limits_{l=0}^{%
+\infty }\left( -1\right) ^{l}\sum\limits_{k=0}^{l}\frac{\left( -1\right) ^{k}(T-r) ^{\left( l-k\right) H}%
}{\left( l-k\right) !}h_{l-k}\left( \frac{B_{T
}^{H}-B_{r }^{H}}{(T-r) ^{H}}\right)\psi_{k}^{(r,T)}(D_{T}^{H,l}F).
\end{equation}
\end{theorem}
The proof is provided in the appendix. To see in a concrete way how to present this convergent series for $\tilde{E}\big[F|\mathcal{F}_{r}^{H}\big]$, we take the following two examples.
\begin{example}
\label{EX0}
Consider the random variable $e^{\sigma B_T^H}$ with $T,~\sigma>0$. One can easily check that it verifies Assumption $(\mathcal A)$. We determine $
\tilde{E}\big[e^{\sigma B_T^H}|\mathcal{F}_{r}^{H}\big]$ for $r\in[0,T]$.
\end{example}
By using (\ref{Iii}),
\begin{equation}
\label{Ex1Eq1}
\tilde{E}\big[e^{\sigma B_T^H}|\mathcal{F}_{r}^{H}\big]=\sum\limits_{l=0}^{%
+\infty }\sum\limits_{k=0}^{l}\frac{\left( -1\right) ^{l+k}(T-r) ^{\left( l-k\right) H}%
}{\left( l-k\right) !}h_{l-k}\left( \frac{B_{T
}^{H}-B_{r }^{H}}{(T-r) ^{H}}\right)\psi_{k}^{(r,T)}(D_{T}^{H,l}e^{\sigma B_T^H}).
\end{equation}
Notice, from (\ref{multiS2}), that for all $l\ge0$ and all $0\le k\le l$,
\begin{equation}
\label{Ex1Eq2}
\psi_{k}^{(r,T)}(D_{T}^{H,l}e^{\sigma B_T^H})=\frac{(T^{2H}-r^{2H})^{k}}{k!}\sigma^{k+l}e^{\sigma B_T^H}.
\end{equation}
It follows by (\ref{Ex1Eq1}), (\ref{Ex1Eq2}) and (\ref{defnHermite}) that
\begin{equation}
\label{Ex1Eq3}
\tilde{E}\big[e^{\sigma B_T^H}|\mathcal{F}_{r}^{H}\big]=e^{\sigma B_T^H}\sum\limits_{l=0}^{%
+\infty }\frac{(-\sigma(T-r)^H)^l}{l!}\sum\limits_{k=0}^{l}{l\choose k}\left(-\frac{\sigma(T^{2H}-r^{2H})}{(T-r)^{H}}\right)^k
h_{l-k}\left( \frac{B_{T
}^{H}-B_{r }^{H}}{(T-r) ^{H}}\right).
\end{equation}
By the following property of Hermite polynomials (due to a Taylor expansion): for all $l\in\mathbb N$, and all $x,y\in\mathbb R$,
$$
\sum\limits_{k=0}^{l}{l\choose k}x^k%
h_{l-k}(y)=h_l(x+y),
$$
one obtains
\begin{equation}
\label{Ex1Eq4}
\sum\limits_{k=0}^{l}{l\choose k}\left(-\frac{\sigma(T^{2H}-r^{2H})}{(T-r)^{H}}\right)^k
h_{l-k}\left( \frac{B_{T
}^{H}-B_{r }^{H}}{(T-r) ^{H}}\right)=h_l\left(\frac{-\sigma(T^{2H}-r^{2H})+B_T^H-B_r^H}{(T-r)^H}\right).
\end{equation}
Finally, it results from (\ref{Ex1Eq3}), (\ref{Ex1Eq4}) and (\ref{interchangetildeE}) that
\begin{eqnarray}
\label{EX01}
\tilde{E}\big[e^{\sigma B_T^H}|\mathcal{F}_{r}^{H}\big] &=&e^{\sigma B_T^H}\sum\limits_{l=0}^{%
+\infty }\frac{(-\sigma(T-r)^H)^l}{l!}h_l\left(\frac{-\sigma(T^{2H}-r^{2H})+B_T^H-B_r^H}{(T-r)^H}\right)\nonumber\\
&=&e^{\sigma B_r^H+\frac{\sigma^2(T^{2H}-r^{2H})}{2}}.
\end{eqnarray}
\begin{example}
\label{Ex1} Consider $F=(B_{t}^H)^{2}B_{T}^H$ for some fixed $T>0$ and $t\in[0,T)$. Below we provide a backward Taylor expansion of $\tilde E[F|\mathcal F_r^H]$ with $r\in[0,T]$.
\end{example}
$F$ is polynomial of fBms, therefore  Assumption $(\mathcal A)$ is verified for all $t\in[0,T)$. By the chain rule (\ref{multiD}), one gets
\begin{eqnarray*} &&\!\!\!\!\!\!\!\!\!\!\!\!D_t^HF=2B_t^HB_T^H+(B_t^H)^2;~D_t^{H,2}F=2B_T^H+4B_t^H;~D_t^{H,3}F=6;~D_t^{H,l}F=0~\mbox{for $l\ge4$};\\
&&\!\!\!\!\!\!\!\!\!\!\!\!D_T^HF=(B_t^H)^2;~D_T^{H,l}F=0~\mbox{for $l\ge2$};~D_t^HD_T^HF=D_T^HD_t^HF=2B_t^H;~D_T^HD_t^{H,2}F=2.
\end{eqnarray*}
Using  Theorem \ref{Theorem 3.2.} and elementary calculus leads to:
if $r\in[0,t]$,
\begin{eqnarray}
\label{Ex10}
\!\!\!\!\! &&\tilde{E}\left[ F|\mathcal{F}_{r}^H\right]=F-(t-r)^Hh_1\left(\frac{B_t^H-B_r^H}{(t-r)^H}\right)D_t^HF+\psi_{1}^{(r,t)}(D_t^HF)\nonumber\\
&&-(T-t)^Hh_1\left(\frac{B_T^H-B_t^H}{(T-t)^H}\right)D_T^HF+\psi_{1}^{(t,T)}(D_T^HF)\nonumber\\
&&+\frac{(t-r)^{2H}}{2!}h_2\left(\frac{B_t^H-B_r^H}{(t-r)^H}\right)D_t^{H,2}F-(t-r)^Hh_1\left(\frac{B_t^H-B_r^H}{(t-r)^H}\right)\psi_{1}^{(r,t)}(D_t^{H,2}F)\nonumber\\
&&+(t-r)^H(T-t)^Hh_1\left(\frac{B_t^H-B_r^H}{(t-r)^H}\right)h_1\left(\frac{B_T^H-B_t^H}{(T-t)^H}\right)D_t^HD_T^HF\nonumber\\
&&-(T-t)^H\psi_1^{(r,t)} \left( h_1\left(\frac{B_T^H-B_t^H}{(T-t)^H}\right)D_T^HD_t^HF \right)\nonumber\\
&&-(t-r)^Hh_1\left(\frac{B_t^H-B_r^H}{(t-r)^H}\right)\psi_1^{(t,T)}(D_t^HD_T^HF)-\frac{(t-r)^{3H}}{3!}h_3\left(\frac{B_t^H-B_r^H}{(t-r)^H}\right)D_t^{H,3}F\nonumber\\
&&-\frac{(t-r)^{2H}(T-t)^H}{2!}h_2\left(\frac{B_t^H-B_r^H}{(t-r)^H}\right)h_1\left(\frac{B_T^H-B_t^H}{(T-t)^H}\right)D_t^{H,2}D_T^HF\nonumber\\
&&+(t-r)^H  (T-t)^H h_1\left(\frac{B_t^H-B_r^H}{(t-r)^H}\right) \psi_1^{(r,t)} \left(h_1\left(\frac{B_T^H-B_t^H}{(T-t)^H}\right) D_t^{H,2}D_T^HF \right) \nonumber \\
&&=(B_r^H)^3+B_r^H\left(T^{2H}+2t^{2H}-3r^{2H}-(T-t)^{2H}\right);
\end{eqnarray}
if $r\in(t,T]$,
\begin{eqnarray}
\label{Ex11}
\tilde{E}\left[ F|\mathcal{F}_{r}^H\right]&=&F-(T-r)^Hh_1\left(\frac{B_T^H-B_r^H}{(T-r)^H}\right)D_T^HF+\psi_{1}^{(r,T)}(D_T^HF)\nonumber\\
&=&B_r^H(B_{t}^H)^{2}+B_{t}^H\left( T^{2H}-r^{2H}-(T-t)^{2H}+(r-t)^{2H}\right).\nonumber\\
\end{eqnarray}
It is interesting to remark that the results from (\ref{EX01}), (\ref{Ex10}) and (\ref{Ex11}) are consistent with the classical conditional expectation with respect to standard Bm. That's because the function $H\longmapsto\tilde E[G(B^H)|\mathcal F^H_t]$ is in fact continuous over $[0,+\infty)$ (see Definition 4.1.2 in \cite{Biagini}). This fact also shows the backward Taylor expansion in Theorem \ref{Theorem 3.2.} naturally extends the one with respect to Bm (see \cite{Schellhorn}).

Note that backward Taylor expansion can be useful when the random variable $F$ is in terms of a discretized fBm trajectory. However, it fails to represent $\tilde E[F|\mathcal F_t^H]$ when $F$ is functional of fBm, such as $F=\int_0^tf(u)\ud B_u^H$. An alternative method is an approximation from its discretization:   if the maximum length in the partition $\max\limits_{i,j\in\{1,\ldots,J\},~i\neq j}|t_i-t_j|$ is small enough, one can roughly approximate $\tilde E[F|\mathcal F_t^H]$ numerically by the series given in (\ref{Ii}).

Our second main result is more general, but under a different sufficient condition for convergence. It does give series representation of $F$ when it is functional of continuous trajectories of fBm. From the fact that $\tilde E[F|\mathcal F_t^H]$ is $\mathcal F_t^H$-measurable, one can imagine that the latter value only depends on the trajectories $\{B_s^H\}_{s\in[0,t]}$, no matter how ill-behaved are the trajectories $\{B_{s}^H\}_{s\in(t,+\infty)}$.  This fact inspires us to originally introduce the "frozen path" operator. More precisely, the series representation given in Theorem \ref{Theorem 3.3.} below can be used to numerically  evaluate a fractional conditional expectation by
following a single typical path backward. It is thus an economical alternative to
Monte Carlo simulation, in case where the fractional Hida Malliavin derivatives are not
expensive to calculate numerically.
\subsection{Series Representation via Exponential Formula}
The "frozen path" operator plays a key role in our second main result: exponential formula of fractional conditional expectations.
 \begin{definition}
 \label{path}
Recall that for $t\ge0$, the action of the element $B^H\chi_{[0,t]}\in L^2(\mathbb P^H)$ on $f\in L^2_{\varphi_H}(\mathbb R)$ is defined as
 $$
 <B^H\chi_{[0,t]},f>:=\int_0^tf(s)\ud B_s^H.
 $$
 For any $F\in L^2(\mathbb P^H)$, let $F=G(B^H)$, where $G$ is the operator defined in Definition \ref{Definition 2.3.}. The "frozen path" operator $\gamma^t:~L^2(\mathbb P^H)\rightarrow L^2(\mathbb P^H)$ is defined as:
 $$
 F(\gamma^t):=G(B^H\chi_{[0,t]}).
 $$
\end{definition}
Note that for $F\in L^2(\mathbb P^H)$, the operator $G$ such that $F=G(B^H)$ always exists, as shown in Theorem 3.1.8 in \cite{Biagini}. Now we explain why the "frozen path" operator is well-defined. Let $\mathcal S$ be the set of polynomial cylindrical random variables of the form
 \begin{equation}
 \label{Fp}
 F=p\left(\int_0^{+\infty}f_1(s)\ud B_s^H,\ldots,\int_0^{+\infty}f_n(s)\ud B_s^H\right),
 \end{equation}
 with $p\in \mathcal P$ (set of polynomials) and $f_i\in L^2_{\varphi_H}(\mathbb R_+)$ for all $i=1,\ldots,n$. The "frozen path" operator $\gamma^t:~\mathcal P\rightarrow L^2(\mathbb P^H)$ is thus well-defined as, for all $t\ge0$,
 $$
 F(\gamma^t)=p\left(\int_0^{t}f_1(s)\ud B_s^H,\ldots,\int_0^{t}f_n(s)\ud B_s^H\right).
 $$
 Here it is obvious that $\gamma^t$ is linear and the function $t\longmapsto F(\gamma^t)$ is continuous on $[0,+\infty)$. Since the set of all random variables of the form in (\ref{Fp}) is dense in $L^2(\mathbb P^H)$ (see e.g. Page 27 in \cite{Nualart} and Page 37 in \cite{Biagini}) and the following proposition shows the "frozen path" operator $\gamma^t$ is continuous from $L^2(\mathbb P^H)$ to $L^2(\mathbb P^H)$, therefore one can naturally extend the domain of "frozen path" operator from $\mathcal P$ to $L^2(\mathbb P^H)$, by preserving the above linearity and continuity. The proof of the key proposition below is given in the appendix.
 \begin{proposition}
\label{prop:closable}.
Let $(F_M)_{M\ge1}$ and $F$ belong to $L^2(\mathbb P^H)$. If $F_M\rightarrow F$ in $L^2(\mathbb P^H)$ as $M\rightarrow+\infty$, then for any $t\ge0$,
$$
F_M(\gamma^t)\xrightarrow[M\rightarrow+\infty]{L^2(\mathbb P^H)}F(\gamma^t).
$$
\end{proposition}
\begin{remark} Assume that $F$ is $\mathcal F_T^H$-measurable, it is clear that $F=G(B^H)=G(B^H\chi_{[0,T]})$ and $F(\gamma^t)=G(B^H\chi_{[0,\min\{T,t\}]})$ is $\mathcal F_t^H$-measurable.
\end{remark} In view of Proposition \ref{prop:closable}, it is not difficult to explicitly compute $F(\gamma^t)$ in most cases.
We provide some examples to show how to obtain "frozen paths" of $L^2(\mathbb P^H)$ random variables:
\begin{enumerate}
\item Denote by $p(B_{t_1}^H,B_{t_2}^H,\ldots,B_{t_n}^H)$ a polynomial of fBm. Set $T\ge\max\{t_1,\ldots,t_n\}$, define
$$
G(B^H)=p\left(\int_0^T\chi_{[0,t_1]}(s)\ud B_s^H,\ldots,\int_0^T\chi_{[0,t_n]}(s)\ud B_s^H\right),
$$
then for all $t\ge0$,
$$
\left(p(B_{t_1}^H,\ldots,B_{t_n}^H)\right)(\gamma^t)=p(B_{\min\{t_1,t\}}^H,\ldots,B_{\min\{t_n,t\}}^H).
$$
\item Let $f\in L^2_{\varphi_H}(\mathbb R_+)$. By using Riemann-Stieltjes integral and Proposition \ref{prop:closable},  for all $t,T\ge0$,
$$
\left(\int_0^Tf(s)\ud B_s^H\right)(\gamma^t)=\int_0^{\min\{t,T\}}f(s)\ud B_s^H.
$$
\item For $F=\int_0^TB^H_{s}\ud s$. One has $
F=G(B^H)=\int_{0}^{T}\int_0^s\ud B_u^H\ud s$, and for all $t\in [0,T]$,
$$
\left(\int_0^TB^H_{s}\ud s\right)(\gamma^t)=\int_0^T\int_0^{\min\{s,t\}}\ud B_u^H\ud s=\int_0^{t}B_s^H\ud s+B_t^H(T-t).
$$
\item Since $\int_0^TB_s^H\ud B_s^H$ is of type FWISI, it is incorrect to define $G(B^H)$ as $
G(B^H)=\int_{0}^T\int_0^s\ud B_u^H\ud B_s^H
$.
However, the It\^o formula shows $\int_0^TB_s^H\ud B_s^H=((B_T^H)^2-T^{2H})/2$, then
$$
\left(\int_0^TB_s^H\ud B_s^H\right)(\gamma^t)=\frac{\big(B_{\min\{T,t\}}^H\big)^2-T^{2H}}{2}.
$$
\item For a general element $F\in L^2(\mathbb P^H)$, $F(\gamma^t)$ can be approximated by some sequence of smooth functions
    $$
    \left\{S_n\left(\int_0^tf_1^{(n)}(s)\ud B_s^H,\ldots,\int_0^tf_n^{(n)}(s)\ud B_s^H\right)\right\}_{n\ge1}
    $$
    in $L^2(\mathbb P^H)$, where $S_n\in C^{\infty}(\mathbb R^n)$ and $f_k^{(n)}\in L^2_{\varphi_H}(\mathbb R)$ for all $n\ge 1$, $1\le k\le n$. This is thanks to Proposition \ref{prop:closable} and the fact that the linear span of smooth functionals of fBms is another dense subset of $L^2(\mathbb P^H)$.
\end{enumerate}
In the illustrations Figure 1 and Figure 2, we show how the "frozen" path operator influences the trajectories of stochastic processes.
\begin{figure}[Ht]
\begin{minipage}[b]{0.45\linewidth}
\centering
\includegraphics[width=\textwidth]{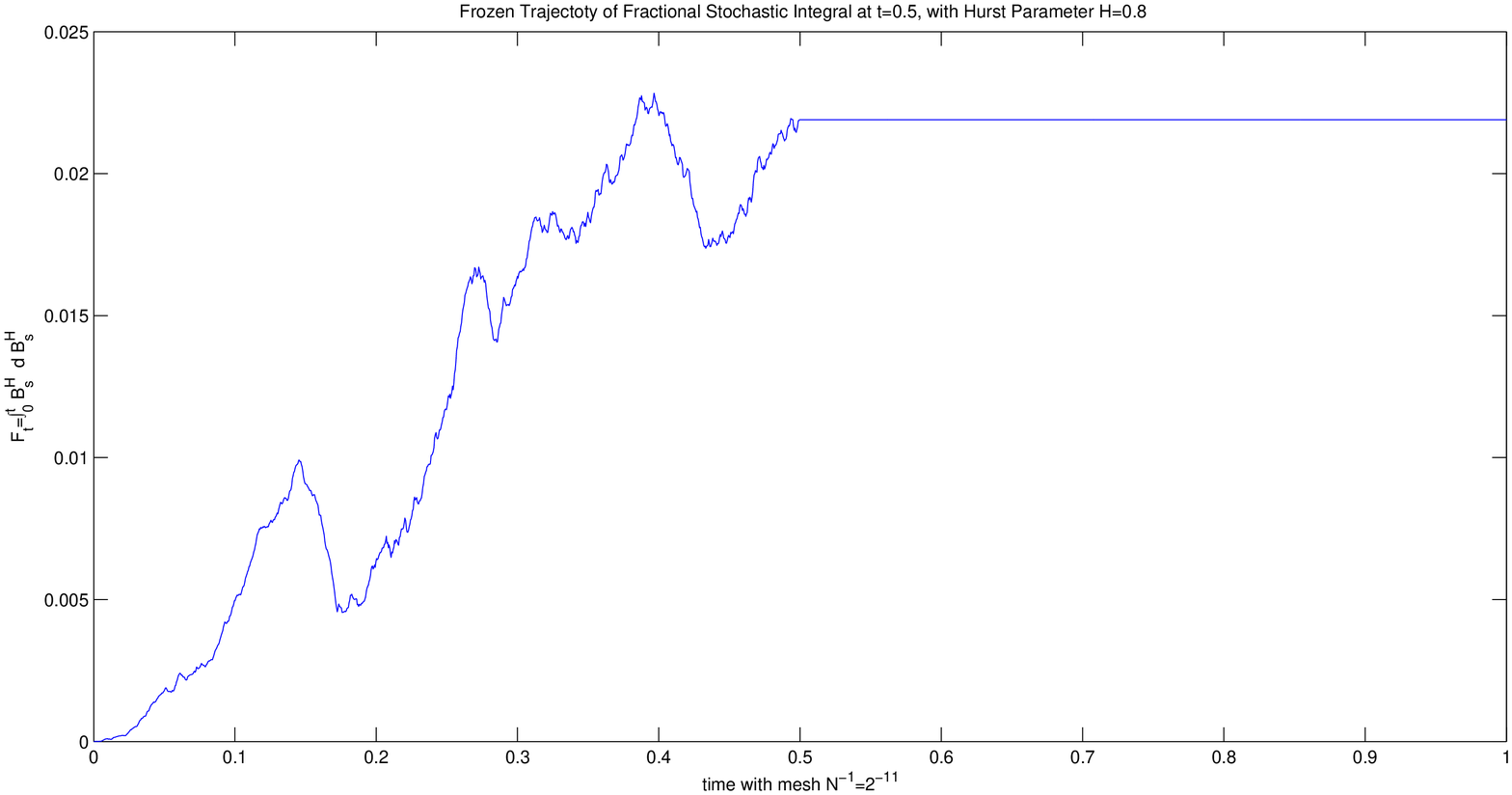}
\caption{Simulation of the frozen path of $\{X_1(s)=\int_0^sB_u^H\ud B_u^H\}_{s\ge0}$, with $H=0.8$, frozen at time $t$=0.5.}
\end{minipage}
\hspace{0.6cm}
\begin{minipage}[b]{0.45\linewidth}
\centering
\includegraphics[width=\textwidth]{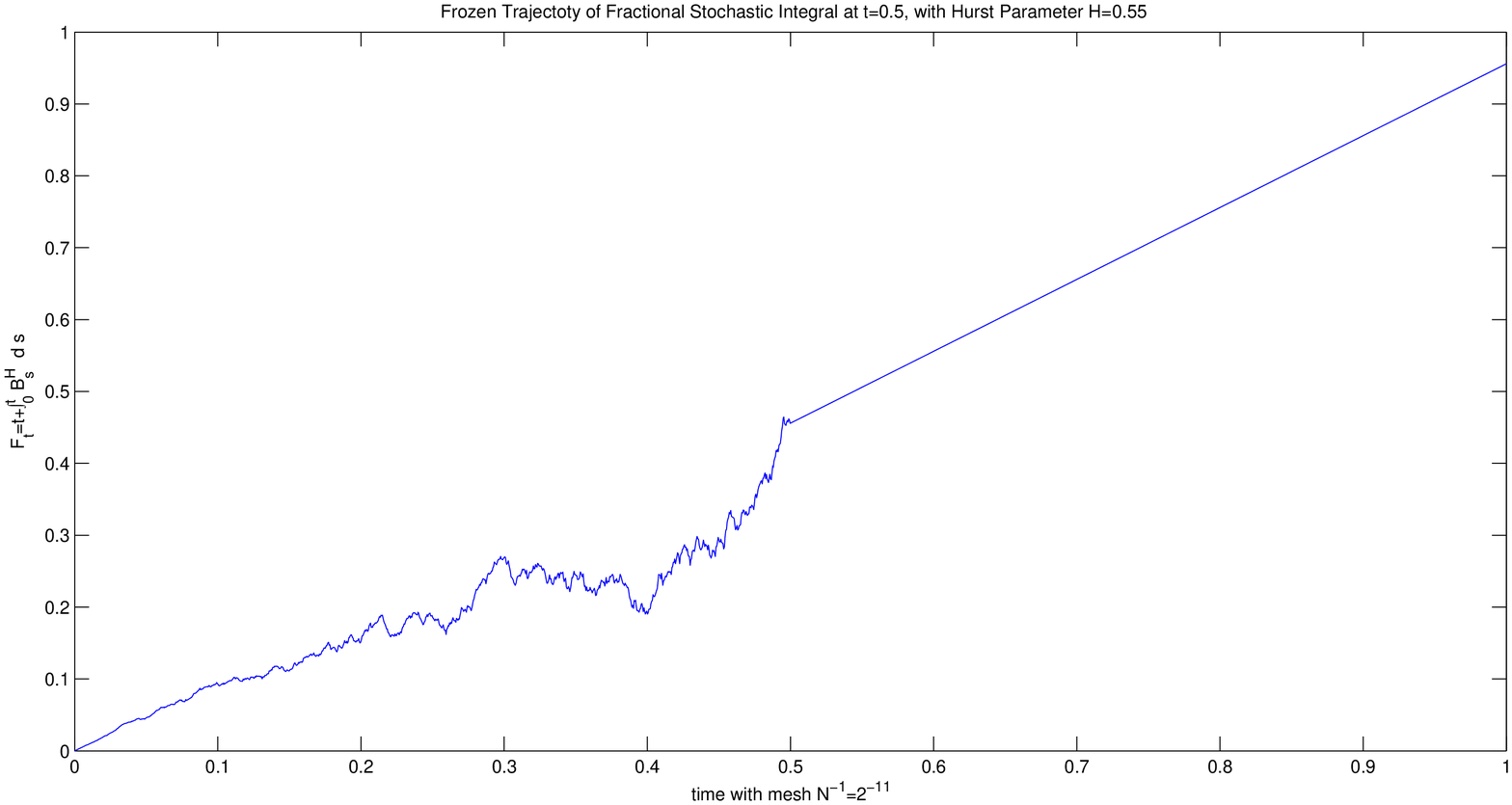}
\caption{Simulation of the frozen path of $\{X_2(s)=s+\int_0^sB_u^H\ud u\}_{s\ge0}$, with $H=0.55$, frozen at time $t$=0.5.}
\end{minipage}
\end{figure}

\noindent\textbf{Assumption ($\mathcal B$):} Let $F\in\mathbb{D}_{\infty,T }^{H}$. Assume that for some $r\in[0, T]$, the following condition holds:
 $$
\sum_{i=1}^{+\infty}\frac{(T^{2H}-r^{2H})^{i}}{2^ii!}\left\|\sup_{u_{2i},\ldots,u_{1}\in[0,T]}\left|\left(D_{u_{2i}}^H\ldots D_{u_1}^HF\right)(\gamma^r)\right|\right\|_{L^2(\mathbb P^H)}<+\infty.
$$
\begin{theorem} [Exponential Formula]
\label{Theorem 3.3.}
For $r\in[0,T]$ and $v\in[r,T]$, define the operator $\mathcal{A}_{v,r}:~L^2(\mathbb P^H)\rightarrow L^2(\mathbb P^H)$ by:
\begin{equation}
\label{Avr}
\mathcal{A}_{v,r}\left( F\right) :=\frac{1}{2}\left(
\int_{0}^TD_{u}^{H}D_{v}^{H}F\varphi_H \left( u,v\right)
\ud u+\int_{0}^{r}D_{u}^{H}D_{v}^{H}F\varphi_H \left( u,v\right) \ud u\right).
\end{equation}
Then under Assumption $(\mathcal B)$, the following series converges in $L^2(\mathbb P^H)$:
\begin{equation}
\label{myequ}
\tilde{E}\left[ F|\mathcal{F}_{r}^H\right] =\sum_{i=0}^{+\infty }\int_{r\leq v_{1}\leq
\ldots\leq v_{i}\leq T}\left( \mathcal{A}_{v_{i},r}...\mathcal{A}_{v_{1},r}F\right)
\left( \gamma^{r}\right) (\ud v)^{\otimes i},
\end{equation}
where by convention, the first item in the series is $F(\gamma^r)$.
\end{theorem}
We emphasize the fact that the integrand of the right-hand side of (\ref{myequ}) must be
evaluated along the path where its driving fBm is frozen at $r$. Below is a quick application of Theorem \ref{Theorem 3.3.}.
\begin{example}
 Consider $F=(B_{t}^H)^{2}B_{T}^H$ for some fixed $T>0$ and $t\in[0,T]$.
 \end{example} By definition, $F$ is $\mathcal F_T^{H}$-measurable and for $r\in[0,t]$, $
F(\gamma ^{r})=(B_{r}^H)^{3}$. Basic computation shows, for any $u,v\in[0,T]$,
$$
\left(D_{u}^HD_{v}^HF\right)(\gamma^r)=2B_{r}^H\left(\chi_{[0,t]}(v)\chi_{[0,t]}(u)+\chi_{[0,t]}(v)\chi_{[0,T]}(u)+\chi_{[0,T]}(v)\chi_{[0,t]}(u)\right).
$$
And then by means of Theorem \ref{Theorem 3.3.}, we get
\begin{eqnarray*}
\tilde{E}\left[ F|\mathcal{F}_{r}^H\right]&=&F\left( \gamma ^{r}\right)+\frac{1%
}{2}\int_{r}^{T}\left( \int_{0}^{T}+\int_0^r\right)\left(D_{u}^{H}D_{v}^{H}F\right)(\gamma^r)\varphi_H \left(
u,v\right) \ud u\ud v\nonumber\\
&=&(B_r^H)^3+B_r^H\left(T^{2H}+2t^{2H}-3r^{2H}-(T-t)^{2H}\right).
\end{eqnarray*}%
This result agrees with (\ref{Ex10}).

Remark that, instead of the backward Taylor expansion, the fractional conditional expectation of the form $F=H(B_{t_1}^H,\ldots,B_{t_J}^H)$ can be also presented via exponential formula, however sometimes it is less obvious to use exponential formula than to use backward Taylor expansion to find the explicit form of $\tilde E[F|\mathcal F_t^H]$. Moreover, backward Taylor expansion is more convenient to use for numerical approximation purpose since it is series of polynomials, while exponential formula contains integrals.

Since $\tilde E[F|\mathcal F_0^H]=E[F]$ (see (\ref{EF=EF}) in the appendix), then one advantage of exponential formula is that it sometimes allows to evaluate expectation of $F$ in an extremely simple way, provided the fractional Malliavin derivatives of $F$ are explicitly known. We motivate this method by giving some real world's applications.
\section{Applications}
\subsection{Fractional Merton Model of Interest Rates}
Suppose that the interest rate follows a particular simple version of the
Merton model: assume $H>1/2$, for $s\in[0,T]$,
\begin{equation*}
r(s)=B_{s}^{H}.
\end{equation*}
The integral $\int_{0}^{T}B_{s}^{H}\ud s$ is defined as a Riemann sum, and the classical result shows
 $$
 \int_{0}^{T}B_{s}^{H}\ud s\sim\mathcal N\left(0,\frac{T^{2H+2}}{2H+2}\right).
 $$
Hence, by using Gaussian moment generating function, the bond price $P(0,T)$ can be directly computed out as
\begin{equation*}
P(0,T)=E\left[ \exp \left(
\int_{0}^{T}B_{s}^{H}\ud s\right) \right]=\exp \left( \frac{T^{2H+2}}{4H+4}\right).
\end{equation*}
We show that this result can be independently obtained by using exponential formula. Note that, applying the exponential formula on $\exp(\int_0^tB_s^H\ud s)$ leads to the same result. To this end, let $F=\exp \big( \int_{0}^{T}B_{s}^{H}\ud s\big)$, it is easy to check that $F$ is $\mathcal F_T^H$-measurable and satisfies Assumption ($\mathcal B$). Observe that
$$
\left(D_{u}^{H}D_{v}^{H}F\right)(\gamma ^{0})=\big(T-u\big)\big(T-v\big),
$$
and then for all integer $i\ge1$,
\begin{eqnarray*}
&&\int_{0\leq v_{1}\leq ...\leq v_{i}\leq T}\left( \mathcal{A}_{v_{i},0}...%
\mathcal{A}_{v_{1},0}F\right) \left( \gamma ^{0}\right) (\ud v)^{\otimes i} \\
&&=\frac{1}{2^{i}}\int_{0\leq v_{1}\leq ...\leq v_{i}\leq
T}\int_{[0,T]^{i}}\prod_{k=1}^{i}(T-u_{k})(T-v_{k})\varphi_H
(u_{k},v_{k})\left( \ud u\right) ^{\otimes i}\left( \ud v\right) ^{\otimes
i}.
\end{eqnarray*}
Since the function $
\int_{[0,T]^{i}}\prod_{k=1}^{i}(T-u_{k})(T-v_{k})\varphi_H
(u_{k},v_{k})\left( \ud u\right) ^{\otimes i}$ is symmetric with respect to $v_1,\ldots,v_i$, one has
\begin{eqnarray*}
&&\int_{0\leq v_{1}\leq ...\leq v_{i}\leq T}\left( \mathcal{A}_{v_{i},0}...%
\mathcal{A}_{v_{1},0}F\right) \left( \gamma ^{0}\right) (\ud v)^{\otimes i} \\
&&=\frac{1}{2^{i}i!}\int_{[0,T]^{2i}}\prod_{k=1}^{i}(T-u_{k})(T-v_{k})%
\varphi_H (u_{k},v_{k})\left( \ud u\right) ^{\otimes i}\left( \ud v\right)
^{\otimes i} \\
&&=\frac{1}{2^{i}i!}\left(\int_{0}^{T}\int_{0}^{T}(T-u)(T-v)\varphi_H (u,v)\ud u\ud v\right)^{i}
\\
&&=\frac{1}{2^{i}i!}\left(\frac{T^{2H+2}}{2H+2}\right)^{i}.
\end{eqnarray*}
Then by applying Theorem \ref{Theorem 3.3.}, we obtain:
\begin{equation*}
P(0,T)=\sum_{i=0}^{+\infty }\frac{1}{2^{i}i!}\left(\frac{T^{2H+2}}{2H+2}%
\right)^{i}=\exp \left( \frac{T^{2H+2}}{4H+4}\right).
\end{equation*}
\subsection{Fractional Cox-Ingersoll-Ross Model of Interest Rates: A Special Case}
Consider a particular case of fractional Cox-Ingersoll-Ross model of interest rates: for $H> 1/2$ and
$s\in [0,T]$, $$r(s)=(B_s^H)^2.$$ Let $F=\exp(-\int_{0}^{T}r(s)\ud s)$, our major goal is to compute $E[F]$. By classical calculus, for any $0\leq u,v,s\leq T$,
\begin{eqnarray*}
D_{v}^{H}F&=&-F\int_{0}^{T} D_{v}^{H}r(s)\ud s=-2F\int_{v}^{T} B_s^H\ud s; \\
D_{u}^{H}D_{v}^{H}F&=&F\int_{0}^{T} D_{u}^{H}r(s)\ud s\int_{0}^{T}D_{v}^{H}r(s)\ud s-F\int_{0}^{T} D_{u}^{H} D_{v}^{H}r(s)\ud s\\
&=&4F\int_{u}^{T}B_s^H\ud s\int_{v}^{T}
B_s^H\ud s-2F\big(T-\max(u,v)\big).
\end{eqnarray*}
Since $r(0)=0$, then imposed to the "frozen path" operator $\gamma ^{0}$, one has,
$$F(\gamma^0)=e^{-\left(\int_0^T(B_s^H)^2\ud s\right)(\gamma^0)}=1;~(D_{u}^{H}D_{v}^{H}F)(\gamma^{0})=-2\big(T-\max(u,v)\big).
$$
Similarly one shows
\begin{eqnarray}
\label{D4F}
&&\left(D_{u_{2}}^{H}D_{v_{2}}^{H}D_{u_{1}}^{H}D_{v_{1}}^{H}F\right)(\gamma^{0})=4\Big(\big(T-\max(u_{1},v_{2})\big)\big(T-\max(u_{2},v_{1})\big)\nonumber \\
&&~~+\big(T-\max(u_{1},u_{2})\big)\big(T-\max(v_{1},v_{2})\big) \nonumber \\
&&~~+\big(T-\max(u_{1},v_{1})\big)\big(T-%
\max(u_{2},v_{2})\big)\Big).
\end{eqnarray}
Then applying Theorem \ref{Theorem 3.3.},
\begin{eqnarray}
\label{item3}
&&E[F]=1+\frac{1}{2}\int_{0}^{T}\int_{0}^{T}
\left(D_{u}^{H}D_{v}^{H}F\right)(\gamma^{0}) \varphi_H(u,v)\ud u\ud v\nonumber \\
&&~~+\frac{1}{4}\int_{[0,T]^3\times[v_1,T]}\left(D_{u_{2}}^{H}D_{v_{2}}^{H}D_{u_{1}}^{H}D_{v_{1}}^{H}F\right)(%
\gamma^{0}) \varphi_H(u_{1},v_{1}) \varphi_H(u_{2},v_{2})\ud v_{2}\ud u_{2}
\ud v_{1}\ud u_{1}\nonumber \\
&&~~+o(T^{4H+2})\nonumber\\
&&=1-\frac{1}{2H+1} T^{2H+1}+\left(\frac{8H^2+18H+5}{4(2H+1)^2(4H+1)}-\frac{\mathcal B(2H+1,2H+2)}{2H+1}\right)T^{4H+2}\nonumber\\
&&~~+o(T^{4H+2}),
\end{eqnarray}
where $\mathcal B$ denotes the beta function, i.e. $\mathcal B(x,y):=\int_0^1t^{x-1}(1-t)^{y-1}\ud t$, for all $x,y>0$. This result provides a numerical approximation of $E[F]$, as $T\rightarrow0$. The computations of the first two items on the right hand-side of the series in (\ref{item3}) are straightforward. However, the third item needs more efforts to obtain. This result is original. We put the technical computation of the third item in the appendix.

Again, note that, although we have assumed $H>1/2$, the expressions of expected values obtained in the previous two examples match the case $H=1/2$ (see \cite{Shreve}), which generalize the results obtained from standard Bm.
\subsection{Series Representation of Characteristic Function of Geometric FBM}
In this example, we show an application of the exponential formula without satisfying Assumption ($\mathcal B$). For $T>0$ and $H>1/2$, let the geometric fBm
$$
X_{T} =e^{\sigma B_{T}^{H}}~\mbox{and}~F=e^{izX_{T}},~\mbox{for $z\in\mathbb C$}.
$$
Recall that $X_T$ is a lognormal random variable. We denote its distribution by $lognormal(0,T^{2H}\sigma^2)$, equivalently, $\log X_T\sim \mathcal N(0,T^{2H}\sigma^2)$. Then $E[F]$ is the characteristic function of some lognormal random variable $%
X_T$ evaluated at $z$. It is known in the literature that there is no convergent Taylor series representation of $E[F]$, due to the fact that the lognormal distribution does not only depend on its probability moments. For example, Holgate \cite{Holgate} showed there is no unique determination of the lognormal distribution only by its moments. Thus, there has been a number of attempts to present the lognormal characteristic functions by divergent power series, which are sufficient for generating moments, e.g. we refer to \cite{Holgate, Barouch,Barakat,Leipnik}. Now we provide a new divergent power series representation of $E[F]$ using Theorem \ref{Theorem 3.3.}. Notice that Assumption $(\mathcal B)$ is no longer satisfied, however the series representation by Theorem \ref{Theorem 3.3.} is still of many interests.

We have, for all $j\in\{1,\ldots,n\}$ and all $0\le u_{j},v_{j}\leq T:$%
\begin{equation*}
D_{u_{j}}^{H}D_{v_{j}}^{H}F=F\sigma ^{2}\big((iz)^{2}e^{2\sigma
B_{T}^{H}}+ize^{\sigma B_{T}^{H}}\big).
\end{equation*}%
By induction,
\begin{equation*}
D_{u_{n}}^{H}D_{v_{n}}^{H}\ldots D_{u_{1}}^{H}D_{v_{1}}^{H}F=F\sigma
^{2n}\sum_{k=0}^{2n}b_{2n,k}(iz)^{k}e^{k\sigma B_{T}^{H}},
\end{equation*}%
where the sequence $\{b_{j,k}\}_{j,k\in\mathbb N}$ satisfies the following recursive formula,
$$b_{j,0}=
\left\{ \begin{array}{ll}
1& \mbox{if $j=0$}\\
0 & \mbox{if $j\neq0$}
\end{array} \right.~\mbox{and for $j\ge k\ge 1$,}~
\left\{\begin{array}{ll}
&b_{j,j} =b_{j,1}=1, \\
&b_{j,k} =kb_{j-1,k}+b_{j-1,k-1}.
\end{array}
\right.
$$
We complete the sequence $\{b_{j,k}\}_{j,k\in\mathbb N}$ by setting $b_{j,k}=0$ for all $0\le j<k$. Thus we claim that it is in fact a sequence of Stirling numbers of the second kind (see for instance \cite{Leipnik} for its definition and properties). We then denote by
$
b_{j,k}=\begin{Bmatrix}
j \\
k%
\end{Bmatrix}$. Therefore by Theorem \ref{Theorem 3.3.} and the fact that $\int_{[0,T]^n}\prod_{l=1}^{n}\varphi_H (u_{l},v_{l}) (\ud u)^{\otimes n} $ is symmetric with respect to $v_1,\ldots,v_n$,
\begin{eqnarray*}
&&\frac{1}{2^{n}}\int_{0\leq v_{1}\leq ...\leq v_{n}\leq
T}%
\int_{[0,T]^{n}}(D_{u_{n}}^{H}D_{v_{n}}^{H}...D_{u_{1}}^{H}D_{v_{1}}^{H}F)(%
\gamma^{0})\prod_{l=1}^n\varphi_H (u_{l},v_{l})(\ud u)^{\otimes
n}(\ud v)^{\otimes n} \\
&&=e^{iz}\frac{1}{2^{n}}\sigma ^{2n}\sum_{k=0}^{2n}%
\begin{Bmatrix}
2n \\
k%
\end{Bmatrix}%
\int_{0\leq v_{1}\leq ...\leq v_{n}\leq
T}\int_{[0,T]^{n}}\prod_{l=1}^{n}\varphi_H (u_{l},v_{l})(\ud u)^{\otimes
n}(\ud v)^{\otimes n} \\
&&=e^{iz}\frac{1}{2^{n}n!}\sigma ^{2n}\sum_{k=0}^{2n}%
\begin{Bmatrix}
2n \\
k%
\end{Bmatrix}%
(iz)^{k}\int_{[0,T]^{n}}\int_{[0,T]^{n}}\prod_{l=1}^{n}\varphi_H
(u_{l},v_{l})(\ud u)^{\otimes n}(\ud v)^{\otimes n} \\
&&=\frac{T^{2nH}}{2^{n}n!}\sigma ^{2n}\sum_{k=0}^{2n}%
\begin{Bmatrix}
2n \\
k%
\end{Bmatrix}%
e^{iz}(iz)^{k},
\end{eqnarray*}
and the power series representation of $E[F]$ is given as
\begin{equation}
\label{CF}
\sum_{n=0}^{+\infty }\sum_{k=0}^{2n}\frac{(\frac{T^{2H}\sigma ^{2}}{%
2})^{n}}{n!}%
\begin{Bmatrix}
2n \\
k%
\end{Bmatrix}%
e^{iz}(iz)^{k}.
\end{equation}
The series (\ref{CF}) is divergent for all $z\in\mathbb C$ due to the large order of the Stirling number of second kind. However, this series representation can be used to evaluate all the moments of $F$, more precisely, by the following relation: for $p\in\mathbb N$,
\begin{equation}
\label{relation}
 E[X_T^p]=(-i)^p\sum_{n=0}^{+\infty }\sum_{k=0}^{2n}\frac{(\frac{T^{2H}\sigma ^{2}}{%
2})^{n}}{n!}i^k
\begin{Bmatrix}
2n \\
k%
\end{Bmatrix}%
\frac{\ud^p(e^{iz}z^{k})}{\ud z^p}\Big|_{z=0}.
\end{equation}
This interesting result provides a new divergent representation of the lognormal distribution's characteristic function using power series. In order to justify this result, we calculate the moments of $X_T$: $ E[X_T^p]$ for $p=1,2,\ldots$ and compare them with the existing results (see e.g. \cite{Leipnik} for Mellin transform of lognormal distribution):
 $$
 E[X_T^p]=\exp\left(p^2 \frac{T^{2H}\sigma^2}{2}\right).
 $$
 Observe that, basic algebraic computation shows for $k\in\mathbb N$,
  $$
  \frac{\ud^p(e^{iz}z^k)}{\ud z^p}\Big|_{z=0}=\sum_{l=0}^pi^{p-l}\frac{p!}{(p-l)!}\delta_l(k),
  $$
  where $\delta$ denotes the Dirac measure by: $\delta_{l}(l)=1$ for $l\in\mathbb Z$, and $\delta_l(k)=0$ for $l,k\in\mathbb Z$ and $l\neq k$.
  Now we take a processing similar to the previous computation and obtain
   \begin{eqnarray*}
  &&\sum_{n=0}^{+\infty }\sum_{k=0}^{2n}\frac{(\frac{T^{2H}\sigma ^{2}}{%
2})^{n}}{n!}%
\begin{Bmatrix}
2n \\
k%
\end{Bmatrix}%
i^{k}\frac{\ud ^p(e^{iz}z^k)}{\ud z^p}\Big|_{z=0}\\
&&=\sum_{n=0}^{+\infty }\sum_{k=0}^{2n}\frac{(\frac{T^{2H}\sigma ^{2}}{%
2})^{n}}{n!}%
\begin{Bmatrix}
2n \\
k%
\end{Bmatrix}%
i^{k}\left(\sum_{l=0}^pi^{p-l}\frac{p!}{(p-l)!}\delta_l(k)\right)\\
&&=i^p\sum_{n=0}^{+\infty }\frac{(\frac{T^{2H}\sigma ^{2}}{%
2})^{n}}{n!}\left(\sum_{l=0}^p\frac{p!}{(p-l)!}\begin{Bmatrix}
2n \\
l%
\end{Bmatrix}\right).
  \end{eqnarray*}
  By the property of generating function of the Stirling number of the second kind, we remark that the following relation holds (see \cite{Graham}) for any $p,n\in\mathbb N$,
  $$
  \sum_{l=0}^p\frac{p!}{(p-l)!}\begin{Bmatrix}
2n \\
l%
\end{Bmatrix}=p^{2n}.
  $$
  It results from (\ref{relation}) that
  $$
  (-i)^p\sum_{n=0}^{+\infty }\sum_{k=0}^{2n}\frac{(\frac{T^{2H}\sigma ^{2}}{%
2})^{n}}{n!}i^k=(-i)^pi^p\sum_{n=0}^{+\infty }\frac{(\frac{T^{2H}\sigma ^{2}}{%
2}p^2)^{n}}{n!}=\exp\left(p^2\frac{T^{2H}\sigma^2}{2}\right).
  $$
  Equivalently, for all $p\in\mathbb N$, $
   E[X_T^p]=\exp(p^2T^{2H}\sigma^2/2)$, which verifies the presentation of characteristic function given in (\ref{CF}) and independently provides another strong demonstration of Theorem \ref{Theorem 3.3.}.
\begin{corollary}
More generally, if $X\sim  lognormal(\mu,\sigma^2)$ with $\mu\in\mathbb R$ and $\sigma>0$, the power series representation of the characteristic function of $X$, $E[e^{izX}]$, can be computed similarly by means of Theorem \ref{Theorem 3.3.} as:
\begin{equation*}
\sum_{n=0}^{+\infty }\sum_{k=0}^{2n}\frac{(\frac{\sigma ^{2}}{%
2})^{n}}{n!}%
\begin{Bmatrix}
2n \\
k%
\end{Bmatrix}%
e^{ize^{\mu}}(ize^{\mu})^{k}.
\end{equation*}
\end{corollary}
\section{Conclusion and Future Work}
We provide two series representations of fractional conditional expectations of functionals of fBm: backward Taylor expansion and exponential formula. We also provide three examples of applications.  The first example is on the fractional Merton model of interest rates, where the series expansion for the bond price can be simplified into a regular exponential, whereas the latter is a well-known result. Our next
two examples are original results. We show how to calculate the first
terms of the series of a bond price in a fractional interest rates model
with time-dependent volatility and are able to approximate the bond price numerically. We also provide a new series representation
of the moment generating function of geometric fBm, which leads to having solved a problem of Fourier transform using power series representation.

So far, our examples only deal with smooth functionals of fBm. As an open problem, it will be interesting, in the future, to develop an approach to represent the "frozen path" of multiple FWISI then give some genius applications of exponential formulas.   Also remark that we have not yet addressed the problem of calculating classical
conditional expectation of a functional of fBm. This work is also quite interesting and meaningful. We note that the methodology developed by
Fourni\'e et al. \cite{Fournie} for Bm
reduces the problem to evaluating two expectations, and is also applicable
to fBm with $H>1/2$. We do not provide explicit
results in that case, and leave this work for future research.
\section{Appendix}
\subsection{Proof of Theorem \ref{Theorem 3.2.}}
One needs the following lemmas to prove Theorem \ref{Theorem 3.2.}. These lemmas, as properties of fractional conditional expectations, generalize the existing results for conditional expectations.
\begin{lemma}
\label{Lemma 4.2.} Let  $%
F\in L^{2}\left( \mathbb{P}^{H}\right) $. For any $t\ge0$, one has
\begin{equation}
\label{tildeE=E}
E\left[ \tilde{E}\left[ F|\mathcal{F}_{t}^{H}\right] \right] =E[F].
\end{equation}
In particular,
\begin{equation}
\tilde{E}\left[ F|\mathcal{F}_{0}^{H}\right] =E\left[ F\right]~\mathbb P^H\mbox{-a.s.}.  \label{EF=EF}
\end{equation}
\end{lemma}
\textbf{Proof of Lemma \ref{Lemma 4.2.}.} For $f\in
L_{\varphi_H}^{2}\left( \mathbb{R_+}\right)$, define the exponential function
\begin{equation}
\label{varepsilonf}
\varepsilon \left( f\right) :=\exp \left( \int_{\mathbb{R_+}}f\left( t\right)
\ud B_{t}^{H}-\frac{1}{2}\left\Vert f\right\Vert _{H,\mathbb R_+}^{2}\right).
\end{equation}
According to Theorem 3.1.4 in \cite{Biagini}, the set of linear span of $%
\varepsilon \left( f\right) $'s is dense in $L^{2}\left( \mathbb{P}^{H}\right) $. Thus it suffices
to prove Lemma \ref{Lemma 4.2.} for $F=\varepsilon \left( f\right)$ then argue with $L^2(\mathbb P^H)$-convergence. Observe that $\varepsilon \left( f\right)$ has the following series representation (we refer to the following statements in \cite{Biagini}: Corollary 3.9.3, Lemma 3.9.2 and Theorem 3.9.7):
$$
\varepsilon \left( f\right) =\sum\limits_{n=0}^{+\infty }\frac{1}{n!}\int_{\mathbb R_+^n}f\left(
s_{1}\right) ...f\left( s_{n}\right) \ud \left(B_{s}^{H}\right)^{\otimes n}.
$$
Thus by definition, for any $t\ge0$,
\begin{equation}
\label{varE}
\tilde{E}\left[ \varepsilon \left( f\right) |\mathcal{F}_{t}^{H}\right]=\sum\limits_{n=0}^{+\infty }\frac{1}{n!}\int_{[0,t]^n}f\left(
s_{1}\right) ...f\left( s_{n}\right)\ud \left(B_{s}^{H}\right)^{\otimes n}\varepsilon\left(f\chi_{[0,t]}\right).
\end{equation}
Recall that, according to Lemma 3.1.3 in \cite{Biagini}, we have, for all $g\in L_{\varphi_H}^2(\mathbb R_+)$,
$
E[ \varepsilon (g)]=1.
$
Since $f$, $f\chi_{[0,t]}\in L_{\varphi_H}^2(\mathbb R_+)$, (\ref{tildeE=E}) holds for $F=\varepsilon(f)$:
 $$
 E\big[\tilde{E}[ \varepsilon (f) |\mathcal{F}_{t}^{H}]\big]=E\big[ \varepsilon (g)\big]=1.
 $$
Now we prove that (\ref{tildeE=E}) holds for all $F\in L^2(\mathbb P^H)$. For arbitrary $F\in L^2(\mathbb P^H)$, there exists a sequence $(f_k^{(n)})_{n\ge1,1\le k\le n}$ in $ L_{\varphi_H}^2(\mathbb R_+)$ and a sequence of real values $(c_k^{(n)})_{n\ge1,1\le k\le n}$, such that
$$
 E\Big|F-\sum_{k=1}^{n}c_k^{(n)}\varepsilon(f_k^{(n)})\Big|^2\xrightarrow[n\rightarrow+\infty]{}0.
$$
Recall the inequality that for any two square integrable random variables $X,Y$,
\begin{equation}
\label{ineqEx}
\big|E[X]-E[Y]\big|\leq E\big|X-Y\big|\leq \sqrt{E\big|X-Y\big|^2}.
\end{equation}
As a consequence,
\begin{equation*}
\Big|E[F]-E\Big[\sum_{k=1}^{n}c_k^{(n)}\varepsilon(f_k^{(n)})\Big]\Big|\leq \sqrt{ E\Big|F-\sum_{k=1}^{n}c_k^{(n)}\varepsilon(f_k^{(n)})\Big|^2}\xrightarrow[n\rightarrow+\infty]{}0.
\end{equation*}
Namely,
\begin{equation}
\label{EF}
E[F]=\lim_{n\rightarrow+\infty}E\Big[\sum_{k=1}^{n}c_k^{(n)}\varepsilon(f_k^{(n)})\Big]=\lim_{n\rightarrow+\infty}\sum_{k=1}^{n}c_k^{(n)}E[\varepsilon(f_k^{(n)})]=\lim_{n\rightarrow+\infty}\sum_{k=1}^{n}c_k^{(n)}.
\end{equation}
The latter limit exists, because $F\in L^2(\mathbb P^H)$. On the other hand, by using the linearity of the fractional conditional expectation and (\ref{boundtildeE}), we obtain
\begin{eqnarray*}
&&E\Big|\tilde{E}\big[ F |\mathcal{F}_{t}^{H}%
\big]-\tilde E\Big[\sum_{k=1}^{n}c_k^{(n)}\varepsilon(f_k^{(n)})\big|\mathcal F_t^H\Big]\Big|^2=E\Big|\tilde{E}\Big[ F-\sum_{k=1}^{n}c_k^{(n)}\varepsilon(f_k^{(n)}) \big|\mathcal{F}_{t}^{H}%
\Big]\Big|^2\\
&&\leq E\Big|F-\sum_{k=1}^{n}c_k^{(n)}\varepsilon(f_k^{(n)})\Big|^2\xrightarrow[n\rightarrow+\infty]{}0.
\end{eqnarray*}
This together with (\ref{ineqEx}) implies
\begin{equation}
\label{EtildeEF}
E\big[\tilde{E}\big[ F |\mathcal{F}_{t}^{H}\big]\big]=\lim_{n\rightarrow+\infty}E\Big[\tilde E\Big[\sum_{k=1}^{n}c_k^{(n)}\varepsilon(f_k^{(n)})\big|\mathcal F_t^H\Big]\Big]=\lim_{n\rightarrow+\infty}\sum_{k=1}^{n}c_k^{(n)}=E[F].
\end{equation}
Therefore (\ref{tildeE=E}) follows from (\ref{EF}) and (\ref{EtildeEF}). In order to obtain (\ref{EF=EF}), it is sufficient to take $t=0$ in (\ref{tildeE=E}) then one notices that $\tilde E[F|\mathcal F_0^H]$ is degenerated to some constant. $\square$

The following lemma extends Proposition 1.2.4 in \cite{Nualart} from Bm to fBm.
\begin{lemma}
\label{lemma:DE}
Let  $%
F\in L^{2}\left( \mathbb{P}^{H}\right)$ be fractional Hida Malliavin differentiable. For any $t,u\ge0$, one has
\begin{equation}
\label{DtildeE=ED}
D_u^H \tilde{E}\left[ F|\mathcal{F}_{t}^{H}\right] =\tilde E\left[D_u^HF|\mathcal F_t^H\right]\chi_{[0,t]}(u).
\end{equation}
\end{lemma}
\textbf{Proof of Lemma \ref{lemma:DE}.} Along the same line as we proved Lemma \ref{Lemma 4.2.}, we first show (\ref{DtildeE=ED}) holds for any $\varepsilon(f)$ defined in (\ref{varepsilonf}), then extend it to linear combinations of $\varepsilon(f)$'s. On the one hand, it follows by (\ref{varE}) and the chain rule that
\begin{eqnarray}
\label{FDE}
D_u^H\tilde E\left[\varepsilon(f)|\mathcal F_t\right]&=&e^{ \int_{0}^tf(s)
\ud B_{s}^{H}-\frac{1}{2}\left\Vert f\chi_{[0,t]}\right\Vert _{H,\mathbb R_+}^{2}}D_u^H\Big(\int_{0}^tf(s)\ud B_s^H\Big)\nonumber\\
&=&\varepsilon(f\chi_{[0,t]})f(u)\chi_{[0,t]}(u).
\end{eqnarray}
On the other hand, again by the chain rule and (\ref{varE}), one gets
\begin{equation}
\label{EFD}
\tilde E\left[D_u^H\varepsilon(f)|\mathcal F_t\right]=\tilde E\left[\varepsilon(f)f(u)\big|\mathcal F_t^H\right]=\varepsilon(f\chi_{[0,t]})f(u).
\end{equation}
Then (\ref{DtildeE=ED}) holds for $F=\epsilon(f)$, thanks to (\ref{FDE}) and (\ref{EFD}). Similar to the proof of (\ref{tildeE=E}), by the denseness of the linear span of $\epsilon(f)$ in $L^2(\mathbb P^H)$, one can show (\ref{DtildeE=ED}) also holds for all fractional Hida Malliavin differentiable $F\in L^2(\mathbb P^H)$. $\square$
\begin{definition}
 Let $F\in \mathbb D_{\infty,T}^H$. For $0\leq r<b\le T$, define the sequence
 $\left\{\left( \int_{r}^{b}\right) ^{H,k}F\ud u\ud v\right\}_{k\in\mathbb N} $ by:
\begin{equation}
\label{multiI1}
\left( \int_{r}^{b}\right) ^{H,0}F\ud u\ud v:=F;
\end{equation}
and for $k\ge1$,
\begin{equation}
\label{multiI2}
\left( \int_{r}^{b}\right) ^{H,k}F\ud u\ud v:=\!\!\!\!\!\!\!\!\!\!\int\limits_{\left[ r\leq v_{k}\leq ...\leq
v_{1}\leq b\right]\times\mathbb{[}0,T\mathbb{]}^{k}}\!\!\!\!\!\!\!\!\!\!D_{u_{k}}^{H}...D_{u_{1}}^{H}F\left(\prod_{i=1}^k\varphi_H \left( u_{i},v_{i}\right)\right) \left( \ud u\right) ^{\otimes k}\left( \ud v\right)
^{\otimes k},
\end{equation}
where $\left( \ud u\right)
^{\otimes k}$ and $\left( \ud v\right)
^{\otimes k}$ are defined as in (\ref{du}).
\end{definition}
\begin{lemma} [Iterated Fractional Integration]
\label{lemma 3.1.}
Let $F\in \mathbb D_{\infty,T}^H$. For any $0\leq r< b\le T$ and any integer $n\ge1$:%
\begin{equation}
\label{intF}
\!\!\!\!\!\!\!\!\!\!\int\limits_{r\leq s_{1}\leq ...\leq s_{n}\leq b}\!\!\!\!\!\!\!\!\!\!F\ud\left( B_{s}^{H}\right)
^{\otimes n}=\sum\limits_{i=0}^{n}\frac{\left( -1\right) ^{i}\left(
b-r\right) ^{\left( n-i\right) H}}{\left( n-i\right) !}h_{n-i}\left( \frac{%
B_{b}^{H}-B_{r}^{H}}{\left( b-r\right) ^{H}}\right) \left(
\int_{r}^{b}\right) ^{H,i}\!\!\!\!F\ud u\ud v,
\end{equation}
where $h_{n-i}\left( x\right) $ is the Hermite polynomial of degree $n-i$ defined in (\ref{Hermite}).
\end{lemma}
\textbf{Proof of Lemma \ref{lemma 3.1.}.}\\
For simplifying notation we define, for any integer $n\ge1$,
$$
a_n:=\int_{r\leq s_{1}\leq ...\leq s_{n}\leq b}F \ud (B_{s}^{H})^{\otimes n}.
$$
We prove Equation (\ref{intF}) by induction. When $n=1$, since $F\in \mathbb D_{\infty,T}^H$, then by (2.33) in \cite{Biagini}, the facts that $h_0\equiv1$, $h_1=Id$ and definitions (\ref{multiI1}), (\ref{multiI2}), we get
\begin{eqnarray*}
a_{1}&=&\int_{r}^{b}F\ud B_{s}^{H}\\
&=&F\int_{r}^{b}\ud B_{s}^{H}-\int_{r}^{b}\int_{0}^{T}D_{u}^{H}F\varphi_H
\left( u,v\right) \ud u\ud v \\
&=&\left( b-r\right) ^{H}h_{1}\left( \frac{B_{b}^{H}-B_{r}^{H}}{\left(
b-r\right) ^{H}}\right) F-\left( \int_{r}^{b}\right) ^{H,1}F\ud u\ud v\\
&=&\sum\limits_{i=0}^{1}\frac{\left( -1\right) ^{i}\left(
b-r\right) ^{\left( 1-i\right) H}}{\left( 1-i\right) !}h_{1-i}\left( \frac{%
B_{b}^{H}-B_{r}^{H}}{\left( b-r\right) ^{H}}\right) \left(
\int_{r}^{b}\right) ^{H,i}F\ud u\ud v.
\end{eqnarray*}
Now assume that (\ref{intF}) holds for $a_k$ with some integer $k\ge1$, then it follows from (\ref{eq}) that
\begin{eqnarray}
\label{inductionhyp}
a_{k}&=&\sum\limits_{i=0}^{k}\frac{\left( -1\right) ^{i}\left(
b-r\right) ^{\left( k-i\right) H}}{\left( k-i\right) !}h_{k-i}\left( \frac{%
B_{b}^{H}-B_{r}^{H}}{\left( b-r\right) ^{H}}\right) \left(
\int_{r}^{b}\right) ^{H,i}F\ud u\ud v\nonumber\\
&=&\sum\limits_{i=0}^{k}\left( -1\right) ^{i}\left(\int_{r\leq s_{1}\leq
s_{2}\leq ...\leq s_{k-i}\leq
b}\ud (B_{s}^{H})^{\otimes (k-i)} \right) \left(
\int_{r}^{b}\right) ^{H,i}F\ud u\ud v.
\end{eqnarray}
Then for $n=k+1$, we write $a_{k+1}$ as
\begin{equation}
\label{ak1}
a_{k+1}=\int_{r}^{b} \left( \int_{s_{1}\leq
s_{2}\leq ...\leq s_{k+1}\leq b}
F\ud (B_{s}^{H})^{\otimes k}\right) \ud B_{s_{1}}^{H}.
\end{equation}
(\ref{ak1}) together with induction hypothesis (\ref{inductionhyp}) yields
\begin{equation*}
a_{k+1}=\sum\limits_{i=0}^{k}(-1)^iI_i,
\end{equation*}
where for $i\in \left\{ 0,1,...,k\right\}$,
\begin{equation}
\label{I}
I_{i}=\int_{r}^{b}\left(\left( \int_{s_{1}}^{b}\right)
^{H,i}F\ud u\ud v\left( \int_{s_{1}\leq s_{2}\leq\ldots\leq s_{k-i+1}\leq
b}\ud (B_{s}^{H})^{\otimes (k-i)}\right)\right)\ud B_{s_{1}}^{H}.
\end{equation}%
Applying Fubini theorem and applying again integration by parts of type for iterated stochastic integrals, we get
\begin{eqnarray*}
I_{i}&=&\!\!\!\!\!\!\!\!\!\!\!\!\!\!\!\!\!\!\!\!\!\!\!\!\int\limits_{\left[ r\leq v_{i}\leq\ldots\leq
v_{1}\leq b\right]\times \mathbb{[}0,T\mathbb{]}^{i} }\!\!\!\!\!\!\!\!\!\!\!\!\int_{r}^{v_i} \Big( D_{u_{i}}^{H}...D_{u_{1}}^{H}F \int_{s_{1}\leq s_{2}\leq ...\leq
s_{k-i+1}\leq b} \ud (B_{s}^{H})^{\otimes (k-i)} \Big) \ud B_{s_{1}}^{H}\\
&&\times\varphi_H \left(
u_{i},v_{i}\right)\ldots \varphi_H \left( u_{1},v_{1}\right)( \ud u) ^{\otimes i}( \ud v)
^{\otimes i}\\
&=&I_{i,1}+I_{i,2},
\end{eqnarray*}%
with
\begin{eqnarray}
\label{I01}
I_{0,1}&=&F\int_{r\leq s_{1}\leq s_{2}\leq ...\leq s_{k+1}\leq
b}\ud (B_{s}^{H})^{\otimes (k+1)};\\
I_{0,2}&=&-\int_r^b\int_0^TD_{u_{1}}^{H}F \int_{v_{1}\leq s_{2}\leq \ldots\leq
s_{k+1}\leq b} \ud (B_{s}^{H})^{\otimes k}\varphi_H \left(
u_{1},v_{1}\right)\ud u_1\ud v_1.\nonumber
\end{eqnarray}
And for $i\in\{1,\ldots,k\}$,
\begin{eqnarray}
\label{Ii1}
I_{i,1}&=&\!\!\!\!\!\!\!\!\!\!\!\!\!\!\!\!\!\!\!\!\!\!\!\!\int\limits_{\left[ r\leq v_{i}\leq \ldots\leq
v_{1}\leq b\right]\times\mathbb{[}0,T\mathbb{]}^{i} }\!\!\!\!\!\!\!\!\!\!\!\! D_{u_{i}}^{H}\ldots D_{u_{1}}^{H}F\int_{r}^{v_i}  \int_{s_{1}\leq s_{2}\leq \ldots\leq
s_{k-i+1}\leq b} \ud (B_{s}^{H})^{\otimes (k-i)} \ud B_{s_{1}}^{H}\nonumber\\
&&\times\varphi_H \left(
u_{i},v_{i}\right)\ldots \varphi_H \left( u_{1},v_{1}\right)( \ud u) ^{\otimes i}( \ud v)
^{\otimes i};\\
\label{Ii2}
I_{i,2} &=&-\!\!\!\!\!\!\!\!\!\!\!\!\!\!\!\!\!\!\!\!\!\!\!\!\int\limits_{\left[r\leq v_{i}\leq \ldots\leq
v_{1}\leq b\right]\times[0,T]^{i} }\!\!\!\!\!\!\!\!\!\!\!\! \Big( \int_{r}^{v_i}\int_{0}^{T}D_{u_{i+1}}^{H}\ldots D_{u_{1}}^{H}F \int_{v_{i+1}\leq s_{2}\leq\ldots\leq
s_{k-i+1}\leq b} \ud (B_{s}^{H})^{\otimes (k-i)}\nonumber\\
&&\times\varphi_H \left( u_{i+1},v_{i+1}\right)
\ud u_{i+1}\ud v_{i+1} \Big) \varphi_H \left(
u_{i},v_{i}\right)\ldots \varphi_H \left( u_{1},v_{1}\right)( \ud u) ^{\otimes i}( \ud v)
^{\otimes i}\nonumber \\
&=&-\!\!\!\!\!\!\!\!\!\!\!\!\!\!\!\!\!\!\!\!\!\!\!\!\int\limits_{\left[ r\leq v_{i+1}\leq \ldots\leq
v_{1}\leq b\right]\times\mathbb{[}0,T\mathbb{]}^{i+1} }\!\!\!\!\!\!\!\!\!\!\!\! D_{u_{i+1}}^{H}\ldots D_{u_{1}}^{H}F \int_{v_{i+1}\leq s_{2}\leq \ldots\leq
s_{k-i+1}\leq b} \ud (B_{s}^{H})^{\otimes (k-i)}\nonumber\\
&&\times\varphi_H \left(
u_{i+1},v_{i+1}\right)\ldots \varphi_H \left( u_{1},v_{1}\right)( \ud u) ^{\otimes (i+1)}( \ud v)
^{\otimes (i+1)}.
\end{eqnarray}
Observe that, on one hand,
\begin{equation}
\label{ak}
a_{k+1}=I_{0,1}+\sum_{i=0}^{k-1}(-1)^i\left( I_{i,2}-I_{i+1,1}\right)+(-1)^kI_{k,2}.
\end{equation}
On the other hand, by the expressions of $I_{i,1}$ in (\ref{Ii1}) and $I_{i,2}$ in (\ref{Ii2}), we have, for $i\in\{0,\ldots,k-1\}$,
\begin{equation}
\label{II}
I_{i,2}-I_{i+1,1}=-\left( \int_{r}^{b}\right)
^{H,i+1}F \ud u\ud v \left( \int_{r\leq s_{1}\leq s_{2}\leq ...\leq s_{k-i}\leq
b}\ud (B_{s}^{H})^{\otimes (k-i)}\right).
\end{equation}
Also notice that the last item in (\ref{ak}) is
\begin{equation}
\label{Ik2}
I_{k,2}=-\left( \int_{r}^{b}\right) ^{H,k+1}F\ud u\ud v.
\end{equation}
It results from (\ref{ak}), (\ref{II}), (\ref{I01}), (\ref{Ik2}) that
\begin{eqnarray*}
a_{k+1}&=&\sum\limits_{i=0}^{k+1}\left( -1\right) ^{i}\left(\int_{r\leq s_{1}\leq
s_{2}\leq ...\leq s_{k+1-i}\leq
b}\ud (B_{s}^{H})^{\otimes (k+1-i)}\right)\left(
\int_{r}^{b}\right) ^{H,i}F\ud u\ud v.
\end{eqnarray*}%
Finally by (\ref{eq}), we obtain
$$
a_{k+1}=\sum\limits_{i=0}^{k+1}\frac{\left( -1\right) ^{i}\left(
b-r\right) ^{\left( k+1-i\right) H}}{\left( k+1-l\right) !}h_{k+1-i}\left( \frac{%
B_{b}^{H}-B_{r}^{H}}{\left( b-r\right) ^{H}}\right) \left(
\int_{r}^{b}\right) ^{H,i}F\ud u\ud v,
$$
which shows (\ref{intF}) holds for $n=k+1$. Therefore Lemma \ref{lemma 3.1.} holds for all integer $n\ge1$. $\square $

\begin{lemma}
\label{LemmaEqual}
Let $F$ be given in Definition \ref{multiintegral}. For $k\in\mathbb N$, $r\in[t_{j-1},t_j)$,
\begin{equation}
\label{equal}
\left(\int_{r}^{t_j}\right) ^{H,k}F\ud u\ud v=\psi_k^{(r,t_j)}(F),
\end{equation}
where $\psi_k^{(r,t_j)}(F)$ is defined in (\ref{multiS1}) and (\ref{multiS2}).
\end{lemma}
\textbf{Proof of Lemma \ref{LemmaEqual}.} When $k=0$, (\ref{equal}) is trivial. When $k\ge1$, by using the fact that $F=H(B_{t_1}^H,\ldots,B_{t_J}^H)$ and (\ref{multiD}), one has, for $u_1,\ldots,u_k\in[0,T]$,
\begin{equation}
\label{multiDu}
D_{u_k}^HD_{u_{k-1}}^H\ldots D_{u_1}^HF=D_{t_1}^{H,\sharp\{u_i:u_i\in[0,t_1],~i=1,\ldots,k\}}\ldots D_{t_J}^{H,\sharp\{u_i:u_i\in(t_{J-1},t_J],~i=1,\ldots,k\}}F,
\end{equation}
where $\sharp\{\cdot\}$ denotes the cardinality of set. It follows from (\ref{multiI2}), (\ref{multiDu}), the multinomial theorem, (\ref{multiS1}) and (\ref{multiS2}) that
\begin{eqnarray}
&&\left( \int_{r}^{t_j}\right) ^{H,k}F\ud u\ud v=\frac{1}{k!}\int_{[r,t_j]^k\times\mathbb{[}0,t_J\mathbb{]}^{k}}D_{u_{k}}^{H}\ldots D_{u_{1}}^{H}F\prod_{i=1}^k\varphi_H \left(
u_{i},v_{i}\right)\left( \ud u\right) ^{\otimes k}\left( \ud v\right)
^{\otimes k}\nonumber\\
&&=\frac{1}{k!}\sum_{\sum_{i=1}^Jq_i=k}{k\choose q_1,q_2,\ldots,q_J}D_{t_1}^{H,q_1}\ldots D_{t_J}^{H,q_J}F\prod_{i=1}^J\left(\int_{r}^{t_j}\int_{t_{i-1}}^{t_i}\varphi_H(u,v)\ud u\ud v\right)^{q_i}\nonumber\\
&&=\psi_k^{(r,t_j)}(F),
\end{eqnarray}
where $
{k\choose q_1,q_2,\ldots,q_J}:=\frac{k!}{q_1!q_2!\ldots q_J!}
$ denotes multinomial coefficient. $\square$

Before introducing the next lemma, to simplify notation, for some proper stochastic process $X:=\{X(s)\}_{s\in \mathbb R^N}$, we indicate its multiple FWISI by
$$
\delta ^{N}(X):=\int_{\mathbb R^N}X(s_1,\ldots,s_N)\ud
\left(B_s^H\right)^{\otimes N}.
$$
 The following technical lemma makes an identification of $E[\delta ^{N}(\mathcal{H}_{N})^{2}]$. It is a particular case of Equation (2.12) in \cite{Nourdin}. More than the latter reference, we provide an explicit form of the norm of tensor product of the Hilbert space $L^2_{\varphi_H}(\mathbb R)$.
\begin{lemma}
\label{deltaN}
For any symmetric stochastic process  $X:=\{X(s)\}_{s\in \mathbb R^N}$ satisfying Assumption $(\mathcal{A})$, we have
\[
E[\delta ^{N}(X)^{2}]=\sum_{i=0}^{N}{N\choose i}^{2}i!E\left\Vert D^{N-i}X \right\Vert _{{L^2_{\varphi_H}(\mathbb R_+)}^{\otimes (2N-i)}}^{2},
\]
where ${L^2_{\varphi_H}(\mathbb R)}^{\otimes (2N-i)}$ denotes the ($2N-i$)th tensor product of $L^2_{\varphi_H}(\mathbb R)$, i.e., $\left\Vert D^{N-i}X \right\Vert _{{L^2_{\varphi_H}(\mathbb R)}^{\otimes (2N-i)}}^{2}$ is given as
\begin{eqnarray*}
&&\left\Vert D^{N-i}X \right\Vert _{{L^2_{\varphi_H}(\mathbb R)}^{\otimes (2N-i)}}^{2}\nonumber \\
&&:=\int_{\mathbb{R}^{4N-2i}}%
\bigg( D_{x_{1},...,x_{N-i}}^{H,N-i}X%
(s_{1},\ldots,s_{N})\prod_{r=1}^{N-i}\varphi_H
(s_{r}^{\prime },x_{r})(\ud x)^{\otimes (N-i)}(\ud s^{\prime})^{\otimes (N-i)}\bigg)\\
&&~~\times\bigg( D_{y_{1},...,y_{N-i}}^{H,N-i}X(s_{1}^{\prime
},\ldots,s_{N}^{\prime
})\prod_{r=1}^{N-i}\varphi_H (s_{r},y_{r})(\ud y)^{\otimes (N-i)}(\ud s)^{\otimes (N-i)}\bigg)\\
&&~~\times\prod_{r=1}^{i}\varphi_H (s_{N-i+r},s_{N-i+r}^{\prime })\ud s_{N}\ldots\ud s_{N-i+1}\ud s'_{N}\ldots\ud s'_{N-i+1}.
\end{eqnarray*}
Here we denote by $D_{x_{1},\ldots,x_{N-i}}^{H,N-i}:=D_{x_1}^H\ldots D_{x_{N-i}}^H;~(\ud x)^{\otimes
(N-i)}:=\ud x_{N-i}\ldots\ud x_1$.
\end{lemma}
Lemma \ref{deltaN} can be obtained without much effort by induction with initial step (which develops the deterministic space $L^2_{\varphi_H}(\mathbb R)$ to space of stochastic processes, see (3.41) in \cite{Biagini}):
\begin{eqnarray*}
&&\left\Vert X \right\Vert _{{L^2_{\varphi_H}(\mathbb R)}}^{2}=\int_{\mathbb R^2}X(s)X(t)\varphi_H(s,t)\ud s\ud t\nonumber\\
&&~~+\int_{\mathbb R^2}\left(\int_{\mathbb R}\varphi_H(s,v)D_v^HX(t)\ud v\right)\left(\int_{\mathbb R}\varphi_H(t,u)D_u^HX(s)\ud u\right)\ud s\ud t.
\end{eqnarray*}
Thus we omit its proof.

\noindent\textbf{Proof of Theorem \ref{Theorem 3.2.}.} Let $r\in[0,T]$ be fixed. Since $\tilde{E}\left[ F|\mathcal{F}_{r}^{H}\right]$ is $\mathcal F_r^H$-measurable, then by the fractional Clark-Hausmann-Ocone formula given in Theorem \ref{Theorem 2.7.},
\begin{equation*}
\tilde{E}\left[ F|\mathcal{F}_{r}^{H}\right] =E\left[ \tilde{E}\left[ F|%
\mathcal{F}_{r}^{H}\right] \right] +\int_{0}^{r}\tilde{E}\left[ D_{s_{1}}^{H}%
\tilde{E}\left[ F|\mathcal{F}_{r}^{H}\right] |\mathcal{F}_{s_{1}}^{H}\right]
\ud B_{s_{1}}^{H}.
\end{equation*}%
Since $r\le T$, by applying Lemma \ref{lemma:DE}, Lemma \ref{Lemma 4.2.} and (\ref{interchangetildeE}),
\begin{eqnarray}
\label{developE}
&&\tilde{E}\left[ F|\mathcal{F}_{r}^{H}\right]\nonumber\\
&&=\left(E\left[ \tilde{E}\left[ F|%
\mathcal{F}_{T}^{H}\right] \right] + \int_{0}^{T}\tilde{E}\left[ D_{s_{1}}^{H}F|\mathcal{F}_{s_{1}}^{H}\right] \ud B_{s_{1}}^{H}\right)-\int_{r}^{T}
\tilde{E}\left[ D_{s_{1}}^{H}F|\mathcal{F}_{s_{1}}^{H}\right] \ud B_{s_{1}}^{H}\nonumber\\
&&=\tilde{E}\left[ F|\mathcal{F}_{T}^{H}\right] -\int_{r}^{T}\tilde{E}\left[
D_{s_{1}}^{H}F|\mathcal{F}_{s_1}^{H}\right] \ud B_{s_{1}}^{H}\nonumber\\
&&=F -\int_{r}^{T}\tilde{E}\left[
D_{s_{1}}^{H}F|\mathcal{F}_{s_1}^{H}\right] \ud B_{s_{1}}^{H},
\end{eqnarray}
because $F$ is $\mathcal F_T^H$-measurable. Repeatedly using the Clark-Hausmann-Ocone formula $n$ times in (\ref{developE}) leads to
\begin{eqnarray}
\label{Clark}
&&\tilde{E}\left[ F|\mathcal{F}_{r}^{H}\right] =F-\int_r^TD_{s_1}^HF\ud B_s^H\nonumber\\
&&~~+\sum_{l=2}^{N-1}(-1)^{l}\int_r^T\int_{s_1}^T\ldots\int_{s_{l-1}}^TD_{s_l}^H\ldots D_{s_1}^HF\ud (B_{s}^H)^{\otimes  l}+R^{(N)},
\end{eqnarray}%
where $R^{(N)}$ is the remainder of the series given by
\begin{equation}
\label{Rndef}
R^{(N)}:=\left( -1\right)
^{N}\int_{r}^{T}\int_{s_{1}}^{T}\ldots\int_{s_{N-1}}^{T}%
\tilde{E}\left[D_{{s_N}}^{H}\ldots D_{s_1}^HF|\mathcal{F}%
_{s_{N}}^{H}\right] \ud (B_{s}^H)^{\otimes N}.
\end{equation}
Next, by the fact that $F=H(B_{t_1}^H,\ldots,B_{t_J}^H)$ and (\ref{multiD}), one has, for $s_1\in(r,T]$,
$$
D_{s_{1}}^{H}F =D_{t_{I_r}}^{H}F\chi_{(r,t_{I_r}]}(s_1)+D_{t_{I_r+1}}^{H}F\chi_{(t_{I_r},t_{I_r+1}]}(s_1)+\ldots+D_{T}^{H}F\chi_{(t_{J-1},t_{J}]}(s_1).
$$
In consequence, for $l\ge1$,
\begin{eqnarray}
\label{developE1}
&&\int_{r\le s_1\le\ldots\le s_l\le T}D_{s_{l}}^{H}\ldots D_{s_1}^HF\ud (B_{s}^H)^{\otimes l}=\!\!\!\!\!\!\!\!\!\!\sum_{q_{I_r}+\ldots+q_J=l}\int_{r\le s_1\le\ldots\le s_l\le T}D_{t_{I_r}}^{H,q_{I_r}}\ldots D_{t_J}^{H,q_J}F\nonumber\\
&&~~\times\chi_{([r,t_{I_r}]^{q_{I_r}}\times(t_{I_r},t_{I_r+1}]^{q_{I_r+1}}\ldots\times(t_{J-1},t_J]^{q_J})}(s_1,\ldots,s_l)\ud (B_{s}^H)^{\otimes l}.
\end{eqnarray}
It results from (\ref{developE1}) and Lemma \ref{lemma 3.1.} that
\begin{eqnarray}
\label{developE2}
&&\int_{r\le s_1\le\ldots\le s_l\le T}D_{s_{l}}^{H}\ldots D_{s_1}^HF\ud (B_{s}^H)^{\otimes l}\nonumber\\
&&=\sum_{q_{I_r}+\ldots+q_{J}=l}\sum_{i_{I_r}=0}^{q_{I_r}}\ldots\sum_{i_J=0}^{q_J}\prod\limits_{k=I_r}^{J}\frac{\left( -1\right) ^{i_k}(\tilde t_{k}-\tilde t_{k-1}) ^{\left( q_{k}-i_k\right) H}%
}{\left( q_{k}-i_k\right) !}\nonumber\\
&&~~\times  h_{q_{I_{r}}-i_{I_{r}}}\left( \frac{B_{\tilde t_{I_{r}}
}^{H}-B_{\tilde t_{I_{r}-1} }^{H}}{(\tilde t_{I_{r}}-\tilde t_{I_{r}-1}) ^{H}} \right) \left(\int_{r}^{t_{I_{r}}}\right)^{H,i_{I_{r}}}\!\!\!\!\!\!\!\ldots
 h_{q_{J}-i_J}\left( \frac{B_{\tilde t_J
}^{H}-B_{\tilde t_{J-1} }^{H}}{(\tilde t_J-\tilde t_{J-1}) ^{H}}\right)\left(\int_{t_{J-1}}^{t_J}\right)^{H,{i_J}} \!\!\!\!\!\!\!\!\!\!\!\! \nonumber\\
&&~~~~D_{t_{I_r}}^{H,q_{I_r}}\!\!\!\!\!\!\ldots D_{t_J}^{H,q_J}F\left(\ud u\ud v\right)^{\otimes (J-I_r+1)}.
\end{eqnarray}
Therefore, by (\ref{developE2}) and Lemma \ref{LemmaEqual}, Theorem \ref{Theorem 3.2.} holds provided that
\begin{equation}
\label{Rnconver}
\|R^{(N)}\|_{L^2(\mathbb P^H)}\xrightarrow[N\rightarrow +\infty ]{ }0.
\end{equation}
Now we show (\ref{Rnconver}) holds under Assumption ($\mathcal A$). Similar to (\ref{developE1}),
\begin{eqnarray*}
R^{(N)}&=&(-1)^N\sum_{q_{I_r}+\ldots+q_J=N}\int_{r\le s_1\le\ldots\le s_N\le T}\tilde E\left[D_{t_{I_r}}^{H,q_{I_r}}\ldots D_{t_J}^{H,q_J}F|\mathcal F_{s_N}^H\right]\\
&&\times\chi_{([r,t_{I_r}]^{q_{I_r}}\times(t_{I_r},t_{I_r+1}]^{q_{I_r+1}}\ldots\times(t_{J-1},t_J]^{q_J})}(s_1,\ldots,s_N)\ud (B_{s}^H)^{\otimes N}.
\end{eqnarray*}
Therefore, to prove (\ref{Rnconver}), it is sufficient to show the sequence below tends to $0$ as $N\rightarrow+\infty$: for any $0\le r<T$ and $q_{I_r},\ldots,q_J\in\{0,\ldots,N\}$,
\begin{eqnarray}
&&\sum_{q_{I_r}+\ldots+q_J=N}\bigg\Vert\int_{r\leq s_1\leq\ldots\leq s_N\leq T}\tilde{%
E}\left[D_{t_{I_r}}^{H,q_{I_r}}\ldots D_{t_J}^{H,q_J}F|\mathcal{F}_{s_{N}}^{H}\right]\nonumber\\
&&~~\times\chi_{([r,t_{I_r}]^{q_{I_r}}\times(t_{I_r},t_{I_r+1}]^{q_{I_r+1}}\ldots\times(t_{J-1},t_J]^{q_J})}(s_1,\ldots,s_N) \ud
\left(B_s^H\right)^{\otimes N}\bigg\Vert_{L^2(\mathbb{P}^H)}\nonumber\\
&&~~\xrightarrow[N\rightarrow+\infty]{}0.  \label{bound}
\end{eqnarray}
To this end we define
\begin{eqnarray*}
&&R_{[r,T]}^{(N)}(q_{I_r},\ldots,q_J):=\int_{r\leq s_1\leq\ldots\leq s_N\leq T}\tilde{%
E}\left[D_{t_{I_r}}^{H,q_{I_r}}\ldots D_{t_J}^{H,q_J}F|\mathcal{F}_{s_{N}}^{H}\right]\nonumber\\
&&~~\times\chi_{([r,t_{I_r}]^{q_{I_r}}\times(t_{I_r},t_{I_r+1}]^{q_{I_r+1}}\ldots\times(t_{J-1},t_J]^{q_J})}(s_1,\ldots,s_N) \ud
\left(B_s^H\right)^{\otimes N}.
\end{eqnarray*}
By a symmetrization of the above integrand,
\begin{equation}
\label{symmeRN}
R_{[r,T]}^{(N)}(q_{I_r},\ldots,q_J)=\frac{(-1)^N}{N!}\int_{[r,T]^N}\mathcal H_N(s_1,\ldots,s_N)\ud
\left(B_s^H\right)^{\otimes N},
\end{equation}
where
\begin{equation}
\label{defH}
\mathcal H_N(s_1,\ldots,s_N):=\sum_{\sigma\in S_N}\tilde{%
E}\left[D_{t_{I_r}}^{H,q_{I_r}}\ldots D_{t_J}^{H,q_J}F|\mathcal{F}_{s_{\sigma(N)}}^{H}\right]\chi_{A}(s_{\sigma(1)},\ldots,s_{\sigma(N)})
\end{equation}
is symmetric with respect to $s_1,\ldots,s_N$ with $S_N$ being the symmetric group on $\{1,2,\ldots,N\}$ and
$$
A:=\left\{(s_{1},\ldots,s_{N})\in[r,t_{I_r}]^{q_{I_r}}\times\ldots\times(t_{J-1},t_J]^{q_J}:~r\le s_{1}\le \ldots\le s_{N} \leq t_{J}\right\}.
$$
According to Lemma \ref{deltaN} (we restrict the processes to nonnegative-time) and Fubini theorem,
\begin{eqnarray}
\label{Edelta11}
&&E\left[\delta ^{N}(\mathcal{H}_{N})^{2}\right]=\sum_{i=0}^{N}{N\choose i}^{2}i!E\left\Vert D^{N-i}\mathcal{H}_N%
\right\Vert _{{L^2_{\varphi_H}(\mathbb R_+)}^{\otimes (2N-i)}}^{2}\nonumber \\
&&=\sum_{i=0}^{N}{N\choose i}^{2}i!\int_{\mathbb{R_+}^{4N-2i}}E\left[ D_{x_{1},...,x_{N-i}}^{H,N-i}%
\mathcal{H}%
_{N}(s_{1},\ldots,s_{N})D_{y_{1},\ldots,y_{N-i}}^{H,N-i}%
\mathcal{H}_{N}(s_{1}^{\prime },\ldots,s_{N}^{\prime })\right]  \nonumber \\
&&~~\times\prod_{r=1}^{N-i}\varphi_H (s_{r}^{\prime },x_{r})\varphi_H
(s_{r},y_{r})\prod_{r=1}^{i}\varphi_H (s_{N-i+r},s_{N-i+r}^{\prime })\nonumber\\
&&~~~~(\ud x)^{\otimes (N-i)}(\ud y)^{\otimes (N-i)}(\ud s)^{\otimes N}(\ud s^{\prime
})^{\otimes N}.
\end{eqnarray}
On the one hand, since $\mathcal H_N(s_1',\ldots,s_N')$ is symmetric, then by the following change of variable in (\ref{Edelta11}): $(s_{N-i+1}',\ldots,s_{N}')\longmapsto (s_{N-i+\sigma''(1)}',\ldots,s_{N-i+\sigma''(i)}')$ for any $\sigma''\in S_i$,  we get for any $i\in\{0,1,\ldots,N\}$,
\begin{eqnarray*}
&&\int_{\mathbb{R_+}^{4N-2i}}E\bigg[ D_{x_{1},...,x_{N-i}}^{H,N-i}%
\mathcal{H}%
_{N}(s_{1},\ldots,s_{N})\nonumber\\
&&~~\times D_{y_{1},\ldots,y_{N-i}}^{H,N-i}%
\mathcal{H}_{N}(s_{1}^{\prime },\ldots,s_{N-i}',s_{N-i+\sigma''(1)}',\ldots,s_{N-i+\sigma''(i)}^{\prime })\bigg]  \nonumber \\
&&~~\times\prod_{r=1}^{N-i}\varphi_H (s_{r}^{\prime },x_{r})\varphi_H
(s_{r},y_{r})\prod_{r=1}^{i}\varphi_H (s_{N-i+r},s_{N-i+\sigma''(r)}^{\prime })\nonumber\\
&&~~~~(\ud x)^{\otimes (N-i)}(\ud y)^{\otimes (N-i)}(\ud s)^{\otimes N}(\ud s^{\prime
})^{\otimes N}\nonumber\\
&&=\int_{\mathbb{R_+}^{4N-2i}}E\left[ D_{x_{1},...,x_{N-i}}^{H,N-i}%
\mathcal{H}%
_{N}(s_{1},\ldots,s_{N})D_{y_{1},\ldots,y_{N-i}}^{H,N-i}%
\mathcal{H}_{N}(s_{1}^{\prime },\ldots,s_{N}^{\prime })\right]  \nonumber \\
&&~~\times \prod_{r=1}^{N-i}\varphi_H (s_{r}^{\prime },x_{r})\varphi_H
(s_{r},y_{r})\prod_{r=1}^{i}\varphi_H (s_{N-i+r},s_{N-i+\sigma''(r)}^{\prime })\nonumber\\
&&~~~~(\ud x)^{\otimes (N-i)}(\ud y)^{\otimes (N-i)}(\ud s)^{\otimes N}(\ud s^{\prime
})^{\otimes N}.
\end{eqnarray*}
It yields
\begin{eqnarray}
\label{Edelta}
&&E\left[\delta ^{N}(\mathcal{H}_{N})^{2}\right]\nonumber\\
&&=\sum_{i=0}^{N}{N\choose i}^{2}\!\!\!\!\!\!\int\limits_{\mathbb{R_+}^{4N-2i}}\!\!\!\!\!\!E\left[ D_{x_{1},...,x_{N-i}}^{H,N-i}%
\mathcal{H}%
_{N}(s_{1},\ldots,s_{N})D_{y_{1},\ldots,y_{N-i}}^{H,N-i}%
\mathcal{H}_{N}(s_{1}^{\prime },\ldots,s_{N}^{\prime })\right]  \nonumber \\
&&~~\times\prod_{r=1}^{N-i}\varphi_H (s_{r}^{\prime },x_{r})\varphi_H
(s_{r},y_{r})\sum_{\sigma''\in S_i}\prod_{r=1}^{i}\varphi_H (s_{N-i+r},s_{N-i+\sigma''(r)}^{\prime })\nonumber\\
&&~~~~(\ud x)^{\otimes (N-i)}(\ud y)^{\otimes (N-i)}(\ud s)^{\otimes N}(\ud s^{\prime
})^{\otimes N}.
\end{eqnarray}
On the other hand, by (\ref{defH}), the linearity of expectation, Cauchy-Schwarz inequality, (\ref{boundtildeE}) and the fact that $F=H(B_{t_1}^H,\ldots,B_{t_J}^H)$, we get
\begin{eqnarray}
\label{ED}
&&E\bigg[ D_{x_{1},\ldots,x_{N-i}}^{H,N-i}\mathcal{H}
_{N}(s_{1},\ldots,s_{N})D_{y_{1},...,y_{N-i}}^{H,N-i}%
\mathcal{H}_{N}(s_{1}^{\prime },\ldots,s_{N}^{\prime })\bigg]\nonumber\\
&&= \sum_{\sigma \in S_{N}}\sum_{\sigma ^{\prime }\in S_{N}}E\bigg[\tilde E\left[D_{x_{1},\ldots,x_{N-i}}^{H,N-i}D_{t_{I_r}}^{H,q_{I_r}}\ldots D_{t_J}^{H,q_J}F|\mathcal{F}_{s_{\sigma
(N)}}^H\right]\nonumber\\
&&~~\times\tilde E\left[D_{y_{1},\ldots,y_{N-i}}^{H,N-i}D_{t_{I_r}}^{H,q_{I_r}}\ldots D_{t_J}^{H,q_J}F|\mathcal{F}_{s_{\sigma
^{\prime }(N)}'}^H\right]\bigg]\nonumber\\
&&~~\times\chi _{[ 0,s_{\sigma (N)}]^{N-i}\times[ 0,s_{\sigma '(N)}^{\prime }]^{N-i}}(x_{1},\ldots,x_{N-i},y_1,\ldots,y_{N-i})  \nonumber \\
&&~~\times\chi_{A^2}(s_{\sigma(1)},\ldots,s_{\sigma(N)},s'_{\sigma'(1)},\ldots,s'_{\sigma'(N)}) \nonumber \\
&&\leq E\bigg[ \sup_{p_{I_r}+\ldots+p_{J}= N-i}\left|
D_{t_{I_r}}^{H,p_{I_r}+q_{I_r}}\ldots D_{t_J}^{H,p_J+q_J}F\right|\bigg]^2 \nonumber\\
&&~~\times\sum_{\sigma \in
S_{N}}\sum_{\sigma ^{\prime }\in S_{N}}\chi _{[ 0,T]^{2N-2i}}(x_{1},\ldots,x_{N-i},y_{1},\ldots,y_{N-i})\nonumber\\
&&~~\times\chi_{A^2}(s_{\sigma(1)},\ldots,s_{\sigma(N)},s'_{\sigma'(1)},\ldots,s'_{\sigma'(N)}).
\end{eqnarray}
Plugging (\ref{ED}) into (\ref{Edelta}), taking $w_i=p_i+q_i$ with $i=I_r,\ldots,J$ and using the fact that $q_{I_r}+\ldots+q_J=N$, we obtain
\begin{eqnarray}
\label{EdeltaB}
&&E[\delta ^{N}(\mathcal{H}_{N})^{2}]\leq\sum_{i=0}^{N}{N\choose i}^{2} E\bigg[\sup_{ w_{I_r}+\ldots+w_J=2N-i}\left|
D_{t_{I_r}}^{H,w_{I_r}}\ldots D_{t_J}^{H,w_J}F\right|\bigg]^2\nonumber\\
&&~~\times\int_{\mathbb{R_+}^{4N-2i}
}\sum_{\sigma \in S_{N}}\sum_{\sigma ^{\prime }\in S_{N}}\chi _{[ 0,T]^{2N-2i}}(x_{1},\ldots,x_{N-i},y_{1},\ldots,y_{N-i})\nonumber\\
&&~~\times\chi_{A^2}(s_{\sigma(1)},\ldots,s_{\sigma(N)},s'_{\sigma'(1)},\ldots,s'_{\sigma'(N)})\nonumber \\
&&~~\times\prod_{r=1}^{N-i}\varphi_H (s_{r}^{\prime },x_{r})\varphi_H
(s_{r},y_{r}) \sum_{\sigma'' \in S_{i}} \prod_{r=1}^{i}\varphi_H (s_{N-i+r},s_{N-i+\sigma'' (r)}^{\prime })\nonumber\\
&&~~~~(\ud x)^{\otimes (N-i)}(\ud y)^{\otimes (N-i)}(\ud s)^{\otimes N}(\ud s^{\prime
})^{\otimes N}.
\end{eqnarray}
Observe that, inside the right-hand side of (\ref{EdeltaB}), for each $(\sigma,\sigma')\in S_N^2$, the integral
\begin{eqnarray*}
&&\!\!\!\!\!\!\int\limits_{\mathbb{R_+}^{2N-2i}
}\!\!\!\!\!\!\chi _{[ 0,T]^{2N-2i}}(x_{1},\ldots,x_{N-i},y_{1},\ldots,y_{N-i})\chi_{A^2}(s_{\sigma(1)},\ldots,s_{\sigma(N)},s'_{\sigma'(1)},\ldots,s'_{\sigma'(N)})\nonumber \\
&&\times\prod_{r=1}^{N-i}\varphi_H (s_{r}^{\prime },x_{r})\varphi_H
(s_{r},y_{r}) \sum_{\sigma'' \in S_{i}} \prod_{r=1}^{i}\varphi_H (s_{N-i+r},s_{N-i+\sigma'' (r)}^{\prime })
(\ud x)^{\otimes (N-i)}(\ud y)^{\otimes (N-i)}
\end{eqnarray*}
returns the same value. Also notice that, for any positive-valued function $f$ and $a\le b$,
\[
\int_{a\leq s_{1} \leq \ldots \leq s_{N} \leq b}\!\!\!\!\!f(s_{1},\ldots,s_{N})(\ud s)^{\otimes{N}}
\leq \int_{{}^{a\leq s_{1} \leq \ldots \leq s_{N-i} \leq b,}_{a\leq s_{N-i+1} \leq \ldots \leq s_{N} \leq b}}f(s_{1},\ldots,s_{N})(\ud s)^{\otimes{N}}.
\]
 Therefore (\ref{EdeltaB}) leads to
\begin{eqnarray}
\label{B1}
&&E[\delta ^{N}(\mathcal{H}_{N})^{2}]\leq\sum_{i=0}^{N}{N\choose i}^{2}(N!)^2 E\left[\sup_{ w_{I_r}+\ldots+w_J=2N-i}\left|
D_{t_{I_r}}^{H,w_{I_r}}\ldots D_{t_J}^{H,w_J}F\right|\right]^2\nonumber\\
&&~~\times\int_{B}\prod_{r=1}^{N-i}\varphi_H (s_{r}^{\prime },x_{r})\varphi_H
(s_{r},y_{r}) \sum_{\sigma'' \in S_{i}} \prod_{r=1}^{i}\varphi_H (s_{N-i+r},s_{N-i+\sigma'' (r)}^{\prime })\nonumber \\
&&~~(\ud x)^{\otimes (N-i)}(\ud y)^{\otimes (N-i)}(\ud s)^{\otimes N}(\ud s^{\prime
})^{\otimes N} \nonumber \\
&&=\sum_{i=0}^{N}{N\choose i}^{2}(N!)^2 E\left[\sup_{ w_{I_r}+\ldots+w_J=2N-i}\left|
D_{t_{I_r}}^{H,w_{I_r}}\ldots D_{t_J}^{H,w_J}F\right|\right]^2\frac{i!}{((N-i)!)^2 (i!)^2}\nonumber\\
&&~~\times\int_{[0,T]^{2N-2i}\times[r,T]^{2N}}\prod_{r=1}^{N-i}\varphi_H (s_{r}^{\prime },x_{r})\varphi_H
(s_{r},y_{r}) \prod_{r=1}^{i}\varphi_H (s_{N-i+r},s_{N-i+r}^{\prime })\nonumber \\
&&~~(\ud x)^{\otimes (N-i)}(\ud y)^{\otimes (N-i)}(\ud s)^{\otimes N}(\ud s^{\prime
})^{\otimes N} \nonumber \\
&&=\sum_{i=0}^{N}{N\choose i}^{2}\frac{(N!)^2i!}{((N-i)!)^2 (i!)^2} E\left[\sup_{ w_{I_r}+\ldots+w_J=2N-i}\left|
D_{t_{I_r}}^{H,w_{I_r}}\ldots D_{t_J}^{H,w_J}F\right|\right]^2 \nonumber \\
&&~~\times\left( \frac{%
T^{2H}-r^{2H}+(T-r)^{2H}}{2}\right) ^{2N-2i}\left( T-r\right) ^{2iH},
\end{eqnarray}
where $B:=\{(x,y,s,s')\in [0,T]^{4N-2i}:~x,y\in [0,T]^{N-i},r\leq s_{1}\leq \ldots\leq s_{N-i}\leq T,r\leq s'_{1}\leq \ldots\leq s'_{N-i}\leq T ,r\leq s_{N-i+1}\leq \ldots\leq s_{N}\leq T,r\leq s'_{N-i+1}\leq \ldots\leq s'_{N}\leq T \}$.  Combining (\ref{B1}) and (\ref{symmeRN}), one gets
\begin{eqnarray}
\label{boundRNT}
&&E\left[R_{[r,T]}^{(N)}(q_{I_r},\ldots,q_J)\right]^2 \leq\sum_{i=0}^{N} E\left[\sup_{ w_{I_r}+\ldots+w_J=2N-i}\left|
D_{t_{I_r}}^{H,w_{I_r}}\ldots D_{t_J}^{H,w_J}F\right|\right]^2\nonumber\\
&&~~\times {N\choose i}%
^{4}\frac{i!}{(N!)^{2}}\left( \frac{T^{2H}-r^{2H}+(T-r)^{2H}}{2}\right)
^{2N-2i}\left( T-r\right) ^{2iH}.
\end{eqnarray}
Now we upper bound $\|R^{{(N)}}\|_{L^2(\mathbb P^H)}$. Observe that
$$
\sharp\left\{(q_{I_r},\ldots,q_J):~q_{I_r}+\ldots+q_J=N\right\}={N+J-I_r \choose J-I_r}\le c\frac{(N+J-1)!}{N!},
$$
where $c>0$ is some constant which does not depend on $r$ nor on $N$.
 And since $I_r\ge1$,
 $$
  E\bigg[\sup_{ \sum_{j=I_r}^Jw_{j}=2N-i}\left|
D_{t_{I_r}}^{H,w_{I_r}}\ldots D_{t_J}^{H,w_J}F\right|\bigg]^2\le  E\left[\sup_{ \sum_{j=1}^Jw_{j}=2N-i}\left|
D_{t_{1}}^{H,w_{1}}\ldots D_{t_J}^{H,w_J}F\right|\right]^2.
 $$
 Then, it follows from (\ref{bound}), (\ref{symmeRN}), (\ref{boundRNT}) and the triangle inequality that
\begin{eqnarray*}
&&\Big\|R^{{(N)}}\Big\|_{L^2(\mathbb P^H)}\le \sum_{q_{I_r}+\ldots+q_J=N}\left\|R_{[r,T]}^{(N)}(q_{I_r},\ldots,q_J)\right\|_{L^2(\mathbb P^H)}\\
&&\le c \frac{(N+J-1)!}{N!}\sum_{i=0}^{N} \left\|\sup_{ \sum_{j=1}^Jw_j=2N-i}\left|
D_{t_{1}}^{H,w_{1}}\ldots D_{t_J}^{H,w_J}F\right|\right\|_{L^2(\mathbb P^H)}\nonumber\\
&&~~\times {N\choose i}%
^{2}\frac{(i!)^{1/2}}{2^{N-i}N!}\left( T^{2H}-r^{2H}+(T-r)^{2H}\right)
^{N-i}\left( T-r\right) ^{iH},
\end{eqnarray*}
which tends to $0$ as $N\rightarrow+\infty$, thanks to Assumption $(\mathcal A)$. $\square$

\subsection{Proof of Proposition \ref{prop:closable}}
The following lemmas are useful to the proof of Proposition \ref{prop:closable}.
\begin{lemma}
\label{lemma uniform}
Let $F\in L^2(\mathbb P^H)$ be $\mathcal F_T^H$-measurable. There exist a sequence of smooth functions $S_n\in C^{\infty}(\mathbb R^n)$ and a sequence of functions $(f_k^{(n)})_{n\ge 1,1\le k\le n}$ in $L^2_{\varphi_H}(\mathbb R_+)$, such that
\begin{eqnarray}
\label{Kn}
&&\sup_{r\in[0,T]}\left\|S_n\left(\int_0^rf_1^{(n)}(s)\ud B_s^H,\ldots,\int_0^rf_n^{(n)}(s)\ud B_s^H\right)- G(B^H\chi_{[0,\min\{r,T\}]})\right\|_{L^2(\mathbb P^H)} \nonumber \\
&&~~\xrightarrow[n\rightarrow+\infty]{}0.
\end{eqnarray}
\end{lemma}
\textbf{Proof of Lemma \ref{lemma uniform}.} Since $F$ is $\mathcal F^H_T$-measurable, there exists $G$ such that $F=G(B^H\chi_{[0,T]})$. Recall that the function $r\longmapsto G(B^H\chi_{[0,\min\{r,T\}]})$ is continuous over $\mathbb R_+$, then by the continuity of fBm and the Weierstrass approximation theorem, there is a sequence of polynomials of fBm $$p_n(B^H_{\min\{r,T/n\}},B^H_{\min\{r,2T/n\}},\ldots,B^H_{\min\{r,(n-1)T/n\}},B_{\min\{r,T\}}^H)$$ such that for all $r\in[0,T]$,
\begin{equation}
\label{approx_Poly}
\sup_{r\in[0,T]}\left\|p_n(B^H_{\min\{r,T/n\}},\ldots,B_{\min\{r,T\}}^H)-G(B^H\chi_{[0,\min\{r,T\}]})\right\|_{L^2(\mathbb P^H)}\xrightarrow[n\rightarrow+\infty]{}0.
\end{equation}
It is easy to see, for each $n\ge1$ and $k\in\{1,\ldots,n\}$ and all $r\in[0,T]$, $B^H_{\min\{r,kT/n\}}$ can be approximated in $L^2(\mathbb P^H)$ as
 \begin{equation}
 \label{approx_point}
\sup_{r\in[0,T]}\left\|\left(m\left(e^{m^{-1}\int_0^r\chi_{[0,kT/n]}\ud B^H_s}-1\right)\right)-B^H_{\min\{r,kT/n\}}\right\|_{L^2(\mathbb P^H)}\xrightarrow[m\rightarrow+\infty]{}0.
 \end{equation}
 Therefore,
  \begin{eqnarray}
  \label{limit3}
  &&\sup_{r\in[0,T]}\left\|p_n\left(\left(n\left(e^{n^{-1}\int_0^r\chi_{[0,kT/n]}\ud B^H_s}-1\right)\right)_{1\le k\le n}\right)-G(B^H\chi_{[0,\min\{r,T\}]})\right\|_{L^2(\mathbb P^H)}\nonumber\\
  &&~~\xrightarrow[n\rightarrow+\infty]{}0.
  \end{eqnarray}
 Observe
 $$
 n\left(e^{n^{-1}\int_0^r\chi_{[0,kT/n]}\ud B^H_s}-1\right)=-n+ne^{\frac{(kT/n)^{2H}}{2n^2}}e^{n^{-1}\int_0^r\chi_{[0,kT/n]}\ud B^H_s-\frac{(kT/n)^{2H}}{2n^2}}.
 $$
  Then there exist two sequences of real values $(a_k^{(n)})_{n\ge1,1\le k\le n}$, $(b_k^{(n)})_{n\ge1,1\le k\le n}$ and a sequence of $L^2_{\varphi_H}(\mathbb R_+)$ functions $f_k^{(n)}=\sum_{i=1}^kb_i^{(n)}\chi_{[0,iT/n]}/n$, so that
$$
\sup_{r\in[0,T]}\left\|\sum_{k=1}^na_k^{(n)}e^{\int_0^tf_k^{(n)}(s)\ud B_s^H-\frac{1}{2}\|f_k^{(n)}\|_{H,\mathbb R_+}^2}- G(B^H\chi_{[0,\min\{r,T\}]})\right\|_{L^2(\mathbb P^H)}\xrightarrow[n\rightarrow+\infty]{}0.
$$
Thus (\ref{Kn}) is obtained when taking $
S_n(x_1,\ldots,x_n)=\sum_{k=1}^na_k^{(n)}e^{x_k-\|f_k^{(n)}\|_{H,\mathbb R_+}^2/2}$. Notice that the exponential functions $S_n$ form a total subset of $C^{\infty}(\mathbb R)$, Lemma \ref{lemma uniform} has been proven. $\square$

 The following preliminary result is a consequence of the fractional Wiener It\^o chaos expansion theorem, it extends Theorem 3.1.8 in \cite{Biagini} and Theorem 2.6 in \cite{Ito}.
 \begin{theorem}
 \label{thm:expansion}
 For $F\in L^2(\mathbb P^H)$, there exists a sequence of real numbers $(c_{\alpha})_{\alpha\in \mathcal J}$, such that for all $r\ge0$,
 \begin{equation}
 \label{expansion1}
 F(\gamma^r)=\sum_{\alpha\in \mathcal J}c_\alpha\tilde{\mathcal H}_{\alpha}(\gamma^r),~\mbox{in $L^2(\mathbb P^H)$},
\end{equation}
 and
 \begin{equation}
 \label{variance1}
 \|F\|_{L^2(\mathbb P^H)}^2=\sum_{\alpha\in \mathcal J}\alpha !c_{\alpha}^2,
 \end{equation}
 where
 \begin{description}
 \item[$-$] $\mathcal J:=\{(\alpha_1,\ldots,\alpha_m)\in\mathbb N^m:~m\ge1\}$ is the set of all finite sequence of nonnegative integers;
 \item[$-$] $\tilde{\mathcal H}_{\alpha}(\gamma^r):=\prod\limits_{i=1}^mh_{\alpha_i}(<B^H\chi_{[0,r]},e_i>)$, with $(e_i)_{i\ge1}$ being the orthonormal basis of $L^2_{\varphi_H}(\mathbb R_+)$ given in (3.10) of \cite{Biagini} and $h_{\alpha_i}$ being the Hermite polynomial defined in (\ref{Hermite});
     \item[$-$] for $\alpha=(\alpha_1,\ldots,\alpha_m)$, $\alpha!:=\alpha_1!\alpha_2!\ldots \alpha_m!$.
 \end{description}
 \end{theorem}
 \textbf{Proof of Theorem \ref{thm:expansion}.} Notice that (\ref{variance1}) is straightforwardly given in Theorem 3.1.8 of \cite{Biagini}, thus we only prove (\ref{expansion1}). Assume that $F$ is $\mathcal F_T^H$-measurable, then for $r> T$, Theorem \ref{thm:expansion} holds thanks to Theorem 3.1.8 in \cite{Biagini}. Now assume $r\in[0,T]$. Set
 $$
 F_n=\sum_{k=1}^na_k^{(n)}\varepsilon(f_k^{(n)}),
 $$
 for $(a_k^{(n)})_{1\le k\le n}\in\mathbb R^n$ and $f_k^{(n)}\in L^2_{\varphi_H}(\mathbb R_+)$. It is shown that (see (3.15) in \cite{Biagini}) for each $f_k^{(n)}$, there exists a sequence of constants $c_{\alpha,k}^{(n)}$ such that
 $$
 \varepsilon(f_k^{(n)})(\gamma^r)=\sum_{\alpha\in\mathcal J}c_{\alpha,k}^{(n)}\tilde{\mathcal H}_{\alpha}(\gamma^r),~\mbox{in $L^2(\mathbb P^H)$}.
 $$
 Therefore Theorem \ref{thm:expansion} holds for $F_n$:
 $$
 F_n(\gamma^r)=\sum_{k=1}^na_k^{(n)}\varepsilon(f_k^{(n)})(\gamma^r)=\sum_{\alpha\in\mathcal J}\left(\sum_{k=1}^na_k^{(n)}c_{\alpha,k}^{(n)}\right)\tilde{\mathcal H}_{\alpha}(\gamma^r).
 $$
 In a general case when $F\in L^2(\mathbb P^H)$ is arbitrary, we apply Lemma \ref{lemma uniform} to  claim that there exists a sequence of functions which is uniformly convergent to $F(\gamma^r)$ in $L^2(\mathbb P^H)$, with respect to $r\in[0,T]$:
 $$
F_n(\gamma^r)=\sum_{k=1}^na_k^{(n)}\varepsilon(f_k^{(n)})(\gamma^r)\xrightarrow[n\rightarrow+\infty]{L^2(\mathbb P^H)}F(\gamma^r).
 $$
 Also observe that, by (3.15) in \cite{Biagini},
 $$
 F_n=\sum_{\alpha\in\mathcal J}\left(\sum_{k=1}^na_k^{(n)}c_{\alpha,k}^{(n)}\right)\tilde{\mathcal H}_{\alpha},
 $$
and by Example 3.1.9 in \cite{Biagini}, the fact that $\lim_{n\rightarrow+\infty}F_n=F$ in $L^2(\mathbb P^H)$, for $\alpha\in\mathcal J$,
 $$
 \sum_{k=1}^na_k^{(n)}c_{\alpha,k}^{(n)}=\frac{1}{\alpha !}E\left[F_n\tilde{\mathcal H}_\alpha\right]\xrightarrow[n\rightarrow+\infty]{}\frac{1}{\alpha !}E\left[F\tilde{\mathcal H}_\alpha\right].
 $$
 Therefore,
 \begin{eqnarray*}
 F(\gamma^r)&=&\lim_{n\rightarrow+\infty}\sum_{\alpha\in\mathcal J}\left(\sum_{k=1}^na_k^{(n)}c_{\alpha,k}^{(n)}\right)\tilde{\mathcal H}_{\alpha}(\gamma^r)\nonumber\\
 &=&\sum_{\alpha\in\mathcal J}\lim_{n\rightarrow+\infty}\left(\sum_{k=1}^na_k^{(n)}c_{\alpha,k}^{(n)}\right)\tilde{\mathcal H}_{\alpha}(\gamma^r)\nonumber\\
 &=&\sum_{\alpha\in\mathcal J}\frac{1}{\alpha !}E\left[F\tilde{\mathcal H}_\alpha\right]\tilde{\mathcal H}_{\alpha}(\gamma^r),~\mbox{in~$L^2(\mathbb P^H)$}.
 \end{eqnarray*}
 It follows that Theorem \ref{thm:expansion} is proven, if we take $c_\alpha=\frac{1}{\alpha !}E[F\tilde{\mathcal H}_\alpha]$. $\square$

Remark that, unlike $(\tilde{\mathcal H}_{\alpha})_{\alpha\in\mathcal J}$ in Theorem 3.1.8 of \cite{Biagini}, $(\tilde{\mathcal H}_{\alpha}(\gamma^r))_{\alpha\in\mathcal J}$ is not necessarily an orthogonal family. In fact, if $F$ is $\mathcal F_T^H$-measurable, then $(\tilde{\mathcal H}_{\alpha}(\gamma^r))_{\alpha\in\mathcal J}$ becomes a family of orthogonal elements in $L^2(\mathbb P^H)$, provided $r\ge T$.

\noindent\textbf{Proof of Proposition \ref{prop:closable}.} Let $(F_M)_{M\ge1}$ and $F$ satisfy $F_M\xrightarrow[M\rightarrow+\infty]{L^2(\mathbb P^H)}F$.
By (\ref{expansion1}), there exists a sequence of constants $(c_{\alpha}^{(M)})_{M\ge1,\alpha\in\mathcal J}$ such that for all $r\ge0$,
\begin{equation}
\label{uniform1}
(F_M-F)(\gamma^r)=\sum_{\alpha\in \mathcal J}c_{\alpha}^{(M)}\tilde{\mathcal H}_{\alpha}(\gamma^r),~\mbox{in $L^2(\mathbb P^H)$}.
\end{equation}
Notice that, by (\ref{variance1}),
\begin{equation*}
\label{uniform2}
\|F_M-F\|_{L^2(\mathbb P^H)}^2=\sum_{\alpha\in \mathcal J}\alpha !(c_{\alpha}^{(M)})^2\xrightarrow[M\rightarrow+\infty]{}0.
\end{equation*}
This together with the fact that $\alpha!>0$ yields $\lim_{M\rightarrow+\infty}c_{\alpha}^{(M)}=0$ for all $\alpha\in\mathcal J$. From (\ref{uniform1}), we see that for any $M\ge1$,
\begin{equation}
\label{uniform4}
\|(F_M-F)(\gamma^r)\|_{L^2(\mathbb P^H)}^2=E\left[\sum_{\alpha\in \mathcal J}\sum_{\beta\in\mathcal J}c_\alpha^{(M)}c_\beta^{(M)}\tilde{\mathcal H}_{\alpha}(\gamma^r)\tilde{\mathcal H}_{\beta}(\gamma^r)\right]<+\infty.
\end{equation}
Since for all $\alpha,\beta\in\cal J$, $
\lim_{M\rightarrow+\infty}|c_\alpha^{(M)}c_\beta^{(M)}\tilde{\mathcal H}_{\alpha}(\gamma^r)\tilde{\mathcal H}_{\beta}(\gamma^r)|=0$. Therefore by (\ref{uniform4}) and the dominated convergence theorem,
$$
\|(F_M-F)(\gamma^r)\|_{L^2(\mathbb P^H)}^2\xrightarrow[M\rightarrow+\infty]{}0.
$$
Proposition \ref{prop:closable} is proven. $\square$
\subsection{Proof of Theorem \ref{Theorem 3.3.}}
Let $F\in\mathbb D_{\infty,T}^H$, since the linear span of $\varepsilon(f\chi_{[0,T]})$ for $f\in L_{\varphi_H}^2(\mathbb R_+)$ is dense in $\mathbb D_{\infty,T}^H$, we first show Theorem \ref{Theorem 3.3.} holds for $F=\varepsilon(f\chi_{[0,T]})$. Applying (\ref{Avr}) to $F=\varepsilon(f\chi_{[0,T]})$ leads to, for any $v\ge0$, $r\in[0,T]$,
\begin{eqnarray}
\label{Avepsilon}
\mathcal{A}_{v,r}\left( \varepsilon(f\chi_{[0,T]})\right)&=&\frac{1}{2}
\left(\int_{0}^{T}+\int_0^r\right)D_{u}^{H}D_{v}^{H}\varepsilon(f\chi_{[0,T]})\varphi_H \left( u,v\right)
\ud u\nonumber\\
&=&\frac{\varepsilon \left( f\chi_{[0,T]}\right)}{2}\left(
\int_{0}^{T}+\int_0^r\right)f(u)f(v)\varphi_H \left( u,v\right)
\ud u.\nonumber\\
\end{eqnarray}
Repeatedly using (\ref{Avepsilon}) implies that for any $v_1,\ldots,v_i\ge0$,
\begin{eqnarray}
\label{Avepsilon1}
&&\left(\mathcal{A}_{v_{i},r}\ldots\mathcal{A}_{v_{1},r}\varepsilon( f\chi_{[0,T]}) \right)\left( \gamma
^{r}\right)\nonumber \\
&&=\frac{\varepsilon \left( f\chi_{[0,T]}\right) \left( \gamma ^{r}\right)}{%
2^{i}}\prod_{k=1}^i\left(\left( \int_{0}^{T}+\int_{0}^{r}\right) f\left(
u\right) f\left( v_{k}\right)\varphi_H \left( u,v_{k}\right)\ud u\right).\nonumber\\
\end{eqnarray}
Observing that $(v_1,\ldots,v_i)\longmapsto \left(\mathcal{A}_{v_{i},r}\ldots\mathcal{A}_{v_{1},r}\left(\varepsilon( f\chi_{[0,T]})\right) \right)\left( \gamma
^{r}\right)$ is symmetric, we obtain from (\ref{Avepsilon1}) that, the following series converges in $L^2(\mathbb P^H)$:
\begin{eqnarray}
\label{avs}
&&\sum_{i=0}^{N }\int_{r\leq v_{1}\leq ...\leq v_{i}\leq T}\left(\mathcal{A}%
_{v_{i},r}...\mathcal{A}_{v_{1},r}\varepsilon\left( f\chi_{[0,T]}\right)\right) \left( \gamma ^{r}\right) (\ud v)^{\otimes i} \nonumber \\
&&=\sum_{i=0}^{N }\frac{1}{i!}\int_{[r,T]^{i}}\left(\mathcal{A}%
_{v_{i},r}...\mathcal{A}_{v_{1},r}\varepsilon\left( f\chi_{[0,T]}\right)\right) \left( \gamma ^{r}\right)(\ud v)^{\otimes i}\nonumber \\
&&=\varepsilon \left( f\chi_{[0,T]}\right) \left( \gamma ^{r}\right) \sum_{i=0}^{N }%
\frac{1}{i!}\left( \frac{1}{2}\int_{r}^{T}\left(
\int_{0}^{T}+\int_{0}^{r}\right) f\left( u\right) f\left( v\right) \varphi_H
\left( u,v\right) \ud u\ud v\right) ^{i}\nonumber\\
&&\xrightarrow[N\rightarrow+\infty]{L^2(\mathbb P^H)}\varepsilon \left( f\chi_{[0,T]}\right) \left( \gamma ^{r}\right)\exp\left(\frac{1}{2}\|f\chi_{[0,T]}\|_{H,\mathbb R_+}^2-\frac{1}{2}\|f\chi_{[0,r]}\|_{H,\mathbb R_+}^2\right).
\end{eqnarray}
Now we determine $\varepsilon \left( f\chi_{[0,T]}\right) \left( \gamma ^{r}\right)$. It follows from (\ref{avs}) and the fact that $\varepsilon(f\chi_{[0,T]})(\gamma^r)
=\exp(\int_0^r f(s)\ud B_s^H-\|f\chi_{[0,T]}\|_{H,\mathbb R_+}^2/2)$ that
\begin{eqnarray}
\label{L2}
&&\sum_{i=0}^{N}\int_{r\leq v_{1}\leq ...\leq v_{i}\leq T}\left(\mathcal{A}%
_{v_{i},r}...\mathcal{A}_{v_{1},r}\left( \varepsilon(f\chi_{[0,T]})\right)\right) \left( \gamma ^{r}\right)(\ud v)^{\otimes i}\nonumber \\
&&\xrightarrow[N\rightarrow+\infty]{L^2(\mathbb P^H)}e^{\int_0^r f(s)\ud B_s^H-\frac{1}{2}\|f\chi_{[0,T]}\|_{H,\mathbb R_+}^2}e^{\frac{1}{2}\|f\chi_{[0,T]}\|_{H,\mathbb R_+}^2-\frac{1}{2}\|f\chi_{[0,r]}\|_{H,\mathbb R_+}^2}\nonumber\\
&&=\varepsilon(f\chi_{[0,r]}).
\end{eqnarray}
By (\ref{varE}) and (\ref{L2}), Theorem \ref{Theorem 3.3.} holds for $F=\varepsilon(f\chi_{[0,T]})$.

We turn to the general case when $F\in\mathbb D_{\infty,T}^{H}$ is arbitrary. Since $r\longmapsto F(\gamma^r)$ is continuous, by Lemma \ref{lemma uniform}, there exists a sequence of coefficients $(a_{i}^{(M)})_{1\le i\le M}$ such that: for all $r\ge0$,
\begin{equation}
\label{uniformFMF}
F_M(\gamma^r):=\sum_{i=1}^{M }a_{i}^{(M)}\varepsilon (f_{i}^{(M)}\chi_{[0,T]})(\gamma^r)\xrightarrow[M\rightarrow+\infty]{L^2(\mathbb P^H)}F(\gamma^r).
\end{equation}
On the one hand, by (\ref{L2}) and the linearity of fractional conditional expectation, one gets
\begin{equation}
\label{conver1}
\sum_{i=0}^{N}\int_{r\leq v_{1}\leq ...\leq v_{i}\leq T}\left(\mathcal{A}%
_{v_{i},r}\ldots\mathcal{A}_{v_{1},r}F_M\right)\left( \gamma ^{r}\right)(\ud v)^{\otimes i}\xrightarrow[N\rightarrow+\infty]{L^2(\mathbb P^H)}\tilde{E}[F_{M}|\mathcal{F}_{r}^{H}].
\end{equation}
The fact that (see (\ref{boundtildeE})) $$
E\left|\tilde{E}[F_{M}-F|\mathcal{F}_{r}^{H}]\right|^2\leq E\left|F_M-F\right|^2\xrightarrow[M\rightarrow+\infty]{}0$$
leads to
\begin{equation}
\label{conver2}
\tilde{E}[F_{M}|\mathcal{F}_{r}^{H}]\xrightarrow[M\rightarrow+\infty]{L^2(\mathbb P^H)}\tilde E[F|\mathcal F_r^H].
\end{equation}
Therefore, from (\ref{conver1}) and (\ref{conver2}) we see the following convergence holds in $L^2(\mathbb P^H)$:
\begin{equation}
\label{cover4}
\lim_{M\rightarrow+\infty}\lim_{N\rightarrow+\infty}\sum_{i=0}^{N}\int_{r\leq v_{1}\leq ...\leq v_{i}\leq T}\left(\mathcal{A}%
_{v_{i},r}\ldots\mathcal{A}_{v_{1},r}F_M\right)\left( \gamma ^{r}\right)(\ud v)^{\otimes i}=\tilde E[F|\mathcal F_r^H].
\end{equation}
On the other hand, one can show that the above $\mathcal{A}%
_{v_{i},r}\ldots\mathcal{A}_{v_{1},r}$ is a closable operator for almost every $v_1,\ldots,v_i\in[r,T]$. In fact, since for almost every $u\ge0$, $D_u^H$ is a closable operator from $L^2(\mathbb P^H)$ to $L^2(\mathbb P^H)$ (see e.g. \cite{Biagini}, Page 38), then by the dominated convergence theorem and (\ref{uniformFMF}), for almost every $v\in [r,T]$,
\begin{eqnarray*}
&&\mathcal{A}_{v,r}F_M =\frac{1}{2}\left(
\int_{0}^{T}D_{u}^{H}D_{v}^{H}F_M\varphi_H \left( u,v\right)
\ud u+\int_{0}^{r}D_{u}^{H}D_{v}^{H}F_M\varphi_H \left( u,v\right) \ud u\right)\nonumber\\
&&\xrightarrow[M\rightarrow+\infty]{L^2(\mathbb P^H)}\frac{1}{2}\left(
\int_{0}^{T}D_{u}^{H}D_{v}^{H}F\varphi_H \left( u,v\right)
\ud u+\int_{0}^{r}D_{u}^{H}D_{v}^{H}F\varphi_H \left( u,v\right) \ud u\right)\nonumber\\
&&=\mathcal A_{v,r}F.
\end{eqnarray*}
By induction, for almost all $v_1,\ldots,v_i\in[r,T]$, $$
\mathcal{A}%
_{v_{i},r}\ldots\mathcal{A}_{v_{1},r}F_{M}\xrightarrow[M\rightarrow+\infty]{L^2(\mathbb P^H)}\mathcal{A}%
_{v_{i},r}\ldots\mathcal{A}_{v_{1},r}F.$$ Further, it follows from Proposition \ref{prop:closable} that, for all $i\ge1$ and almost all $v_1,\ldots,v_i\in[r,T]$,
\begin{equation}
\label{conver3}
\left(\mathcal{A}%
_{v_{i},r}\ldots\mathcal{A}_{v_{1},r}F_{M}\right)( \gamma ^{r})\xrightarrow[M\rightarrow+\infty]{L^2(\mathbb P^H)}\left(\mathcal{A}%
_{v_{i},r}\ldots\mathcal{A}_{v_{1},r}F\right)( \gamma ^{r}).
\end{equation}
Now to show the exponential formula in Theorem \ref{Theorem 3.3.} converges, it suffices to demonstrate that the series
\begin{equation}
\label{convergeA}
\sum_{i=1}^{+\infty}\int_{r\leq v_{1}\leq ...\leq v_{i}\leq T}\left(\mathcal{A}%
_{v_{i},r}...\mathcal{A}_{v_{1},r}F\right)( \gamma ^{r})(\ud v)^{\otimes i}
\end{equation}
is convergent in $L^2(\mathbb P^H)$. To this end one observes, from the binomial theorem, that for $i\ge1$,
\begin{eqnarray*}
&&\left(\mathcal A_{v_i,r}\ldots\mathcal A_{v_1,r}F\right)(\gamma^r)\\
&&=\frac{1}{2^i}\sum_{k=0}^i{i\choose k}\int\limits_{[0,T]^k\times[0,r]^{i-k}}\left(D_{u_i}^HD_{v_i}^H\ldots D_{u_1}^HD_{v_1}^HF\right)(\gamma^r)\prod_{j=1}^i\varphi_H(u_j,v_j)(\ud u)^{\otimes i}\\
&&\le \left(\sup_{u_{2i},\ldots,u_1\in[0,T]}\left|\left(D_{u_{2i}}^H\ldots D_{u_1}^HF\right)(\gamma^r)\right|\right)\frac{1}{2^i}\sum_{k=0}^i{i\choose k}\!\!\!\!\!\!\!\int\limits_{[0,T]^k\times[0,r]^{i-k}}\!\!\!\prod_{j=1}^i\varphi_H(u_j,v_j)(\ud u)^{\otimes i}.
\end{eqnarray*}
It yields
\begin{eqnarray}
\label{BoundSum2}
&&\sum_{i=1}^{N}\int_{r\leq v_{1}\leq ...\leq v_{i}\leq T}\left(\mathcal A_{v_i,r}\ldots\mathcal A_{v_1,r}F\right)(\gamma^r)(\ud v)^{\otimes i}\nonumber\\
&&\le \sum_{i=1}^{N}\left(\sup_{u_{2i},\ldots,u_1\in[0,T]}\left|\left(D_{u_{2i}}^H\ldots D_{u_1}^HF\right)(\gamma^r)\right|\right)\frac{1}{2^i}\sum_{k=0}^i{i\choose k}\nonumber\\
&&~~\times\frac{1}{i!}\int\limits_{[r,T]^i\times[0,T]^k\times[0,r]^{i-k}}\varphi_H(u_i,v_i)\ldots\varphi_H(u_1,v_1)(\ud u)^{\otimes i}(\ud v)^{\otimes i}\nonumber\\
&&=\sum_{i=1}^{N}\frac{1}{2^ii!}\sum_{k=0}^i{i\choose k}\left(\int_r^T\int_0^T\varphi_H(u,v)\ud u\ud v\right)^k
\left(\int_r^T\int_0^r\varphi_H(u,v)\ud u\ud v\right)^{i-k}\nonumber\\
&&~~\times\sup_{u_{2i},\ldots,u_1\in[0,T]}\left|\left(D_{u_{2i}}^H\ldots D_{u_1}^HF\right)(\gamma^r)\right|\nonumber\\
&&=\sum_{i=1}^{N}\frac{(T^{2H}-r^{2H})^i}{2^ii!}\left(\sup_{u_{2i},\ldots,u_1\in[0,T]}\left|\left(D_{u_{2i}}^H\ldots D_{u_1}^HF\right)(\gamma^r)\right|\right).
\end{eqnarray}
Next by applying the triangle inequality to (\ref{BoundSum2}), one obtains
\begin{eqnarray*}
&&\left\Vert\sum_{i=1}^{N}\int_{r\leq v_{1}\leq ...\leq v_{i}\leq T}\left(\mathcal A_{v_i,r}\ldots\mathcal A_{v_1,r}F\right)(\gamma^r)(\ud v)^{\otimes i}\right\Vert_{L^2(\mathbb{P}^H)}\\
&&\le\sum_{i=1}^{N}\frac{(T^{2H}-r^{2H})^i}{2^ii!}
\left\Vert\sup_{u_{2i},\ldots,u_1\in[0,T]}\left|\left(D_{u_{2i}}^H\ldots D_{u_1}^HF\right)(\gamma^r)\right|\right\Vert_{L^2(\mathbb{P}^H)}.
\end{eqnarray*}
The above series is convergent as $N\rightarrow+\infty$, thanks to Assumption ($\mathcal B$).
Finally,  it follows from (\ref{cover4}), (\ref{convergeA}), the triangle inequality and the monotone convergence theorem that the following convergence holds in $L^2(\mathbb P^H)$:
\begin{eqnarray*}
&&\sum_{i=1}^{+\infty}\int_{r\leq v_{1}\leq ...\leq v_{i}\leq T}\left(\mathcal{A}%
_{v_{i},r}\ldots\mathcal{A}_{v_{1},r}F\right)\left( \gamma ^{r}\right)(\ud v)^{\otimes i}\\
&&=\lim_{N\rightarrow+\infty}\lim_{M\rightarrow+\infty}\sum_{i=0}^{N}\int_{r\leq v_{1}\leq ...\leq v_{i}\leq T}\left(\mathcal{A}%
_{v_{i},r}\ldots\mathcal{A}_{v_{1},r}F_{M}\right)\left( \gamma ^{r}\right)(\ud v)^{\otimes i}\\
&&=\lim_{M\rightarrow+\infty}\lim_{N\rightarrow+\infty}\sum_{i=0}^{N}\int_{r\leq v_{1}\leq ...\leq v_{i}\leq T}\left(\mathcal{A}%
_{v_{i},r}\ldots\mathcal{A}_{v_{1},r}F_{M}\right)\left( \gamma ^{r}\right)(\ud v)^{\otimes i}\\
&&=\tilde E[F|\mathcal F_r^H].
\end{eqnarray*}
 Therefore Theorem \ref{Theorem 3.3.} has been proven.
$\square$
\subsection{Computation of (\ref{item3})}
Since the first two items of the expansion in (\ref{item3}) are relatively easy to obtain, in this part we only show that
\begin{eqnarray}
\label{computeItem3}
&&\frac{1}{4}\int_{[0,T]^3\times[v_1,T]}\left(D_{u_{2}}^{H}D_{v_{2}}^{H}D_{u_{1}}^{H}D_{v_{1}}^{H}F\right)(%
\gamma^{0}) \varphi_H(u_{1},v_{1}) \varphi_H(u_{2},v_{2})\ud v_{2}\ud u_{2}
\ud v_{1}\ud u_{1}\nonumber\\
&&=\left(\frac{8H^2+18H+5}{4(2H+1)^2(4H+1)}-\frac{\mathcal B(2H+1,2H+2)}{2H+1}\right)T^{4H+2}.
\end{eqnarray}
\textbf{Proof of (\ref{computeItem3}).} First, from (\ref{D4F}), we observe  that integrating
$$
\left(D_{u_{2}}^{H}D_{v_{2}}^{H}D_{u_{1}}^{H}D_{v_{1}}^{H}F\right)(%
\gamma^{0}) \varphi_H(u_{1},v_{1}) \varphi_H(u_{2},v_{2})
$$
over the following domains gives the same value: $\{(u_1,u_2,v_1,v_2)\in[0,T]^4:~v_2\ge v_1,~u_1\ge v_1,~u_2\ge v_2\}$, $\{(u_1,u_2,v_1,v_2)\in[0,T]^4:~v_2\ge v_1,~u_1\ge v_1,~u_2< v_2\}$, $\{(u_1,u_2,v_1,v_2)\in[0,T]^4:~v_2\ge v_1,~u_1< v_1,~u_2\ge v_2\}$ and  $\{(u_1,u_2,v_1,v_2)\in[0,T]^4:~v_2\ge v_1,~u_1<v_1,~u_2< v_2\}$. As a consequence, the left-hand side of (\ref{computeItem3}) is equal to
\begin{eqnarray}
\label{computeItem4}
&&\frac{1}{4}\int_{[0,T]^3\times[v_1,T]}\left(D_{u_{2}}^{H}D_{v_{2}}^{H}D_{u_{1}}^{H}D_{v_{1}}^{H}F\right)(%
\gamma^{0}) \varphi_H(u_{1},v_{1}) \varphi_H(u_{2},v_{2})\ud v_{2}\ud u_{2}
\ud v_{1}\ud u_{1}\nonumber\\
&&=\!\!\!\!\!\!\!\!\!\!\!\!\int\limits_{[0,T]^4\cap[v_2\ge v_1,u_1\le v_1,u_2\le v_2]}\!\!\!\!\!\!\!\!\!\!\!\!\left(D_{u_{2}}^{H}D_{v_{2}}^{H}D_{u_{1}}^{H}D_{v_{1}}^{H}F\right)(%
\gamma^{0}) \varphi_H(u_{1},v_{1}) \varphi_H(u_{2},v_{2})\ud v_{2}\ud u_{2}
\ud v_{1}\ud u_{1}\nonumber\\
&&=\mathcal I_1+\mathcal I_2+\mathcal I_3,
\end{eqnarray}
where
\begin{eqnarray}
\label{computeI1term}
&&\mathcal I_1=\!\!\!\!\!\!\!\!\!\!\!\!\int\limits_{0\le u_2\le u_1\le v_1\le v_2\le T}\!\!\!\!\!\!\!\!\!\!\!\!\left(D_{u_{2}}^{H}D_{v_{2}}^{H}D_{u_{1}}^{H}D_{v_{1}}^{H}F\right)(%
\gamma^{0}) \varphi_H(u_{1},v_{1}) \varphi_H(u_{2},v_{2})\ud v_{2}\ud u_{2}
\ud v_{1}\ud u_{1}\nonumber\\
&&=4(H(2H-1))^2\nonumber\\
&&~\times\!\!\!\!\!\!\!\!\!\!\!\!\int\limits_{0\le u_2\le u_1\le v_1\le v_2\le T}\!\!\!\!\!\!\!\!\!\!\!\!(v_1-u_1)^{2H-2}(v_2-u_2)^{2H-2}\big[(T-u_1)(T-v_2)+2(T-v_1)(T-v_2)\big]\nonumber\\
&&~~\ud v_{2}\ud u_{2}
\ud v_{1}\ud u_{1}\nonumber\\
&&=\frac{(2H-1)(H+2)}{(2H+1)(4H+2)(4H-1)}T^{4H+2};
\end{eqnarray}
\begin{eqnarray}
\label{computeI2term}
&&\mathcal I_2=\!\!\!\!\!\!\!\!\!\!\!\!\int\limits_{0\le u_1\le u_2\le v_1\le v_2\le T}\!\!\!\!\!\!\!\!\!\!\!\!\left(D_{u_{2}}^{H}D_{v_{2}}^{H}D_{u_{1}}^{H}D_{v_{1}}^{H}F\right)(%
\gamma^{0}) \varphi_H(u_{1},v_{1}) \varphi_H(u_{2},v_{2})\ud v_{2}\ud u_{2}
\ud v_{1}\ud u_{1}\nonumber\\
&&=4(H(2H-1))^2\nonumber\\
&&~\times\!\!\!\!\!\!\!\!\!\!\!\!\int\limits_{0\le u_1\le u_2\le v_1\le v_2\le T}\!\!\!\!\!\!\!\!\!\!\!\!(v_1-u_1)^{2H-2}(v_2-u_2)^{2H-2}\big[(T-u_2)(T-v_2)+2(T-v_1)(T-v_2)\big]\nonumber\\
&&~~\ud v_{2}\ud u_{2}
\ud v_{1}\ud u_{1}\nonumber\\
&&=\left(\frac{8H^2+14H-1}{4(4H+1)(2H+1)(4H-1)}-\frac{6H+5}{2H+1}\mathcal B(2H+1,2H+2)\right)T^{4H+2};\nonumber\\
\end{eqnarray}
and
\begin{eqnarray}
\label{computeI3term}
&&\mathcal I_3=\!\!\!\!\!\!\!\!\!\!\!\!\int\limits_{0\le u_1\le v_1\le u_2\le v_2\le T}\!\!\!\!\!\!\!\!\!\!\!\!\left(D_{u_{2}}^{H}D_{v_{2}}^{H}D_{u_{1}}^{H}D_{v_{1}}^{H}F\right)(%
\gamma^{0}) \varphi_H(u_{1},v_{1}) \varphi_H(u_{2},v_{2})\ud v_{2}\ud u_{2}
\ud v_{1}\ud u_{1}\nonumber\\
&&=4(H(2H-1))^2\nonumber\\
&&~~\times\!\!\!\!\!\!\!\!\!\!\!\!\int\limits_{0\le u_1\le v_1\le u_2\le v_2\le T}\!\!\!\!\!\!\!\!\!\!\!\!(v_1-u_1)^{2H-2}(v_2-u_2)^{2H-2}\big[2(T-u_2)(T-v_2)+(T-v_1)(T-v_2)\big]\nonumber\\
&&~~\ud v_{2}\ud u_{2}
\ud v_{1}\ud u_{1}\nonumber\\
&&=\frac{6H+4}{2H+1}\mathcal B(2H+1,2H+2)T^{4H+2}.
\end{eqnarray}
Finally, (\ref{computeItem3}) follows from (\ref{computeItem4}), (\ref{computeI1term}), (\ref{computeI2term}) and (\ref{computeI3term}). $\square$



\end{document}